\documentclass[11pt]{amsart}
\usepackage{amsmath,amssymb}
\usepackage{mathrsfs}
\usepackage{hyperref}

\newtheorem{mthm}{Theorem}

\newtheorem{thm}{Theorem}[section]
\newtheorem{cor}[thm]{Corollary}
\newtheorem{prop}[thm]{Proposition}

\newtheorem{lem}[thm]{Lemma}

\newtheorem{obs}[thm]{Observation}

\theoremstyle{definition}
\newtheorem{defn}[thm]{Definition}
\newtheorem{defns}[thm]{Definitions}
\newtheorem{exmp}[thm]{Example}

\newtheorem{cons}[thm]{Construction}

\theoremstyle{remark}
\newtheorem{rem}[thm]{Remark}

\def\NN{\mathbf{N}}
\def\ZZ{\mathbf{Z}}
\def\QQ{\mathbf{Q}}
\def\RR{\mathbf{R}}

\def\X{\mathscr{X}}

\newcommand{\overto}[1]{{\buildrel{#1}\over\longrightarrow}}
\newcommand{\overfrom}[1]{{\buildrel{#1}\over\longleftarrow}}

\newcommand{\ol}[1]{\overline{#1}}

\newcommand{\setdef}[2]{ \left\{ {#1}\ :\ {#2} \right\} }

\newcommand{\pa}{\partial}   

\newcommand{\ro}{\varrho}

\def\beq{\begin{equation}}
\def\eeq{\end{equation}}
\def\GL{{\bf\rm GL}}
\def\PGL{{\bf\rm PGL}}
\def\SL{{\bf\rm SL}}
\def\PSL{{\bf\rm PSL}}
\def\SO{{\bf\rm SO}}

\def\supp{{\rm supp}}

\def\Span{{\rm Span}}

\def\Ker{{\rm Ker}}

\def\mod{{/\!\!/}}
\def\pprime{{\prime\prime}}

\def\Id{{\rm Id}}

\def\sA{{\mathcal{A}}}

\def\sH{{\mathcal{H}}}

\def\Prob{{\sf Prob}}
\def\Meas{{\sf Meas}}
\def\Aut{{\rm Aut}\,}

\def\Aff{{\rm Aff}}

\def\Diff{{\rm Diff}}
\def\Isom{{\rm Isom}}
\def\Homeo{{\rm Homeo}}
\def\Stab{{\rm Stab}}
\def\Fix{{\rm Fix}}
\def\Circle{{\rm S}^{1}} 
\def\Sphere{{\rm S}^{2}} 
\def\Sect{{\Sigma}} 
\def\lcsc{{l.c.s.c. }} 

\title[Superrigidity, Weyl Groups, and $\Homeo_{+}(\Circle)$]{Superrigidity, Weyl groups, \\
and actions on the Circle}
\author[Uri Bader]{Uri Bader}
\address{University of Chicago, Chicago, USA}
\email{bader@math.uchicago.edu}

\author[Alex Furman]{Alex Furman}
\address{University of Illinois at Chicago, Chicago, USA}
\email{furman@math.uic.edu}
\thanks{A.F. is supported in part by the NSF grant DMS-0094245.}

\author[Ali Shaker]{Ali Shaker}
\address{University of Illinois at Chicago, Chicago, USA}
\email{shaker@math.uic.edu}

\begin{document}

\begin{abstract}
We propose a new approach to superrigidity phenomena 
and implement it for lattice representations and measurable cocycles 
with $\Homeo_{+}(\Circle)$ as the target group. 
We are motivated by Ghys' theorem stating that
any representation $\ro:\Gamma\to\Homeo_{+}(\Circle)$ of an 
irreducible lattice $\Gamma$ in a semi-simple  real Lie group $G$
of higher rank, either has a finite orbit or, up to a semi-conjugacy,  
extends to $G$ which acts through an epimorphism $G\to\PSL_{2}(\RR)$. 

Our approach, based on the study of abstract boundary theory
and, specifically, on the notion of a generalized Weyl group, allows:
(A) to prove a similar superrigidity result for irreducible
lattices in products $G=G_{1}\times\cdots G_{n}$ of $n\ge 2$ 
general locally compact groups, 
(B) to give a new (shorter) proof of Ghys' theorem,
(C) to establish a commensurator superrigidity for general locally compact groups,
(D) to prove first superrigidity theorems for $\tilde{A}_{2}$ groups.
This approach generalizes to the setting of measurable circle bundles;
in this context we prove cocycle versions of (A), (B) and (D).
This is the first part of a broader project of studying superrigidity via generalized Weyl groups.
\end{abstract}

\maketitle

\begin{center}
{\it Dedicated to Ali\\ (A.F. and U.B.)}
\end{center}

\section{Introduction and Statement of the Main Results}

Celebrated Superrigidity theorem of Margulis 
is one of the central results on lattices in semi-simple Lie groups; 
it gives a complete description of the linear representations of higher rank lattices 
and allows to classify such lattices via Margulis' Arithmeticity theorem
(\cite{Margulis:book:91}). 
Zimmer's Cocycle Superrigidity is another remarkably powerful result which became 
an indispensable tool in the study of  group actions on manifolds and in 
Ergodic theory of actions of large groups (cf. \cite{Zimmer:book:84}).
In recent years the Superrigidity phenomenon has motivated and inspired a lot of further
research on other "higher rank" groups, and the study of representations into target 
groups other than linear ones, e.g. isometry groups of CAT(-1) spaces, 
Hilbert spaces, CAT(0) spaces, $\Homeo(\Circle)$ or $\Diff(\Circle)$. 

Our general objective is to develop a new and unified approach to many of these 
superrigidity phenomena; in the present paper this approach is implemented to 
obtain new superrigidity results (as well as to reprove known ones) for representations 
and cocycles valued in $\Homeo(\Circle)$ -- the group of homeomorphisms of the circle. 
The approach is built around a concept of a \emph{generalized Weyl group} --
which is a certain, typically finite, group $W_{G,B}$ associated to a given locally compact 
group $G$ and (a choice of)  a $G$-\emph{boundary} in the sense of Burger-Monod, i.a.
a measure space with a $G$-action which is amenable and doubly ergodic
with (unitary) coefficients.
For semi-simple $G$ and the natural boundary $B=G/P$ 
the generalized Weyl group agrees with the classical one
(e.g. it is $S_{n}$ for $G=\SL_{n}(k)$ where $k$ is a local field).
Amenable groups $G$ admit a trivial boundary with the corresponding Weyl group being trivial;
if $G$ is an algebraic group of rank-one, or is a general Gromov hyperbolic group
(or a non-elementary subgroup of such) then typically $W_{G,B}\cong \ZZ/2\ZZ$.
However $W_{G,B}$ is more complicated in many situations that we associate 
with ``higher rank'', namely: 
\begin{itemize}
\item[(i)] 
semi-simple algebraic $k$-groups $G$ with ${\rm rank}_{k}(G)\ge 2$; 
\item[(ii)]
products $G=G_{1}\times\cdots\times G_{n}$ of $n\ge 2$ of non-amenable
factors -- in this case the Weyl group $W_{G,B}$ contains $\prod_{i=1}^{n} W_{G_{i},B_{i}}$ 
which contains $(\ZZ/2\ZZ)^{n}$; 
\item[(iii)]
groups acting simply transitively on vertices of $\tilde{A}_{2}$ buildings, in this case 
there exists a boundary with Weyl group containing $S_{3}$.
\end{itemize}

The philosophy is that any complexity of $W_{G,B}$ imposes certain amount of rigidity on
actions of $G$, on representations/actions of lattices $\Gamma<G$, 
and on measurable cocycles of $G$ or $\Gamma<G$.
In the context of actions on the circle, i.e. representations/cocycles into 
$\Homeo(\Circle)$, we shall  see  that rigidity phenomena appear 
whenever $W_{G,B}$ is \emph{not cyclic}.

\bigskip

Actions on the circle, or rather lack of such, for lattices in higher rank Lie groups
have recently attracted the attention of a number of authors.
Most notable are the works of Witte \cite{Witte:circle:94}, 
Ghys \cite{Ghys:lattices:99},  
Burger - Monod \cite{BurgerMonod:JEMS:99}, \cite{BurgerMonod:GAFA:02},
Farb - Shalen \cite{FarbShalen:biggroups:99},
Witte - Zimmer \cite{WitteZimmer:circlebundle:01},
Navas \cite{Navas:Tcircle:02}, \cite{Navas:new1phenomena:05}, \cite{Navas:cocycles:06},
Lifschitz - Morris (Witte) \cite{MorrisLifschitz:Qrank1:04}.
Some of these results require smoothness assumptions, and some
restrict the class of lattices considered. 
In the framework of superrigidity results, which apply to all lattices in a given group $G$,
and assume only continuous action on the circle the best result is due to \'Etienne Ghys:
\begin{thm}[{ \cite[Th\'eor\`eme 3.1]{Ghys:lattices:99}}]
\label{T:Ghys}
Let $G$ be a semi-simple connected real Lie group with finite center 
and ${\rm rank}_{\RR}(G)\ge 2$,
let $\Gamma<G$ be an irreducible lattice, 
and $\ro:\Gamma\to\Homeo_{+}(\Circle)$ a homomorphism. 

Then either $\ro(\Gamma)$ has a finite orbit, 
or $G$ has an epimorphism $p:G\to \PSL_{2}(\RR)$
and the original $\Gamma$-action $\ro$ is semi-conjugate to 
the action
\[
	\ro_{0}:\Gamma< G\overto{p} \PSL_{2}(\RR)<\Homeo_{+}(\Circle)
\]
by fractional linear transformations.
\end{thm}
Let us clarify some of the terms used above. 
Hereafter all the results about actions on the circle are stated for 
$\Homeo_{+}(\Circle)$ -- the group of \emph{orientation preserving}
homeomorphisms of the circle, which has index two in $\Homeo(\Circle)$.
This is done only to simplify the formulations.
We now recall the notion of (semi)-conjugacy and of elementary actions.
\begin{defns}
\label{D:semi-conj}
For two actions $\ro, \ro_{0}:\Gamma\to\Homeo_{+}(\Circle)$ 
we say that $\ro$ is \emph{semi-conjugate} to $\ro_{0}$ 
if there exists a continuous monotonic surjection
$f:\Circle\to\Circle$ so that $\ro_{0}(\gamma)\circ f=f\circ \ro(\gamma)$ 
for all $\gamma\in\Gamma$. 
If $f$ can be chosen to be a homeomorphism, the actions $\ro$ and
$\ro_{0}$ are called \emph{conjugate}. 

An action $\ro:\Gamma\to\Homeo_{+}(\Circle)$
is \emph{elementary} if it has a finite orbit, or is semi-conjugate into a group of rotations.
Equivalently, $\ro(\Gamma)$ has an invariant probability measure on the circle.
\end{defns}
Note that if $H^{1}(\Gamma,\RR)$ is trivial, as is the case for higher rank lattices,
then elementary $\Gamma$-actions are precisely those which have a finite orbit. 
(See subsection~\ref{SS:elementary} below for further discussion).

In the above result of Ghys, one can distinguish between a lattice $\Gamma$ in 
a \emph{simple} Lie group $G$ of higher rank, and the case of an \emph{irreducible} 
lattice $\Gamma$ in a \emph{semi}-simple group 
$G=\prod_{i=1}^{n} G_{i}$ which has $n\ge 2$ non-compact factors.
Non-elementary actions appear only in the latter case, and only
if one of the simple factors is locally isomorphic to $\PSL_{2}(\RR)$.
In the following result we generalize this case to a setting of an \emph{irreducible}
lattice $\Gamma$ in a product $G=\prod_{i=1}^{n} G_{i}$ of $n\ge 2$ \emph{general} \lcsc groups.
In this context \emph{irreducibility} means that $\pi_{i}(\Gamma)$ is dense in $G_{i}$ 
for each $i$, where $\pi_{i}:G\to G_{i}$ denotes the coordinate projection.
\begin{mthm}[Superrigidity for Lattices in General Products]
\label{T:product}
Let $\Gamma$ be an irreducible lattice in a product 
$G=G_{1}\times \cdots\times G_{n}$ of $n\ge 2$  \lcsc groups $G_{i}$, 
which do not have infinite index open normal subgroups.
Let $\ro:\Gamma\to \Homeo_{+}(\Circle)$ be a homomorphism with $\ro(\Gamma)$ 
being non-elementary.

Then one of the factors $G_{i}$ of $G$ admits an epimorphism 
$p:G_{i}\to\PSL_{2}(\RR)$ so that the action $\ro:\Gamma\to\Homeo_{+}(\Circle)$ 
is semi-conjugate to the action 
\[
	\ro_{0}:\Gamma< G\overto{\pi_{i}} G_{i}\overto{p} \PSL_{2}(\RR)<\Homeo_{+}(\Circle)
\]
by fractional linear transformations.
\end{mthm}
%

%
This result has an ``Arithmeticity'' Corollary, when combined with Shalom's \cite[Theorem 0.5]{Shalom:Inven:00}:
\begin{cor}[Arithmeticity]
\label{C:arithm}
Let $\Gamma<G=G_{1}\times \cdots\times G_{n}$ be an irreducible lattice and
$\ro:\Gamma\to\Homeo_{+}(\Circle)$ be a non-elementary action as in Theorem~\ref{T:product}.
Assume that each $G_{i}$ is a simple as a topological group, and that $G/\Gamma$ is
compact.

Then each $G_{i}$ is isomorphic to some $\PSL_{2}(k_{i})$ where $k_{i}$ is
some local fields of zero characteristic, and $k_{j}=\RR$ for some $1\le j\le n$,
$\Gamma$ is an $S$-arithmetic lattice in $G=\prod_{i=1}^{n} \PSL_{2}(k_{i})$ 
and $\ro$ is semi-conjugate to the natural $\Gamma$-action by fractional linear transformations
via the projection to the $j$th factor:
\[
	\Gamma< G\overto{\pi_{j}} \PSL_{2}(\RR)<\Homeo_{+}(\Circle).
\]
\end{cor}
\begin{rem}
The idea that irreducible lattices in products of $n\ge 2$ groups 
should share many rigidity properties with (irreducible) lattices
in higher rank (semi)-simple Lie groups, can be traced back at least to
Benrstein - Kazhdan \cite{BernsteinKazhdan:70} and is present in the works of 
Margulis in the seventies.
This idea is central in the study of lattices in products of automorphism groups of trees 
by Burger - Mozes 
\cite{BurgerMozes:LatProdTrees:00}, 
\cite{BurgerMozes:LocGlob:00}, \cite{BurgerMozes:Simple:97}, see also
\cite{BurgerMozesZimmer:arithtree} and references therein.
Starting from Shalom's work \cite{Shalom:Inven:00} many
rigidity properties of irreducible lattices in products of completely general 
\lcsc groups were obtained;  in this direction see also: Monod - Shalom \cite{MonodShalom:CRAS:03}, 
\cite{MonodShalom:OE:05}, Bader - Shalom \cite{BaderShalom:Factor:06}, 
Monod \cite{Monod:splittingCRAS:05}, \cite{Monod:arithmeticity:05}, \cite{Monod:splittingJAMS:06},
Bader - Furman - Gelander - Monod \cite{BFGM:FTLp:05},
Gelander - Karlsson - Margulis \cite{GelanderKarlssonMargulis}.
The latter preprint also contains a previous result of Margulis on 
abstract commensurator superrigidity \cite{Margulis:commen:94}.
So Theorem~\ref{T:product} is  a natural analogue of Ghys' theorem.

In fact, results similar to Theorem~\ref{T:product} were recently proved by A. Navas in
\cite[Theorems B, C, D]{Navas:new1phenomena:05} using Shalom's splitting 
theorem for affine isometric actions on Hilbert spaces \cite{Shalom:Inven:00}. 
However, in order to apply these techniques it was necessary to impose additional assumptions 
on $\ro$ and on $\Gamma$ in $G$, namely that $\ro(\Gamma)$ is in $\Diff^{1+\tau}(\Circle)$ 
with $\tau>1/2$\footnote{Substituting Shalom's splitting theorem by its $L^{p}$ analogue from \cite{BFGM:FTLp:05} with large $p$,
Navas' arguments can be applied to any $\tau>0$ (cf. \cite{Navas:cocycles:06}).}
and $\Gamma<G$ is cocompact, or satisfies slightly weaker  conditions
as in \cite{Shalom:Inven:00}.
\end{rem}

Let us point out, however, that our approach is different from the 
variety of ideas and the techniques which appear in the above mentioned works on
lattices in general products; it also differs from Ghys' original proof of Theorem \ref{T:Ghys},
and from a later extension to arbitrary local fields obtained by Witte-Zimmer 
\cite{WitteZimmer:circlebundle:01}.
In fact, we can give  an alternative (and, to our taste, a rather easy) proof of theorem~\ref{T:Ghys}.
Observe that Theorem~\ref{T:product} covers irreducible lattices in all semi-simple groups with 
more than one non-compact factor.
Indeed, given such a group $G$, the theorem applies to the universal covering group $\tilde{G}$,
and the resulting map $p\circ \pi_i$ factors through $G$, as $\tilde{G}$ is a central extension of $G$, 
and $\PSL_2(\RR)$ is center free.
It remains to deal with the case of lattices in simple groups.
\begin{mthm}[after Ghys, Witte-Zimmer]
\label{T:higherrank}
Let $k$ be a local field, $\mathbf{G}$ be a $k$-conected simple algebraic group defined over $k$,
${\rm rank}_{k}(\mathbf{G})\ge 2$. Let $\Gamma$ be a lattice 
in the \lcsc group $G=\mathbf{G}(k)$, and $\ro:\Gamma\overto{}\Homeo_{+}(\Circle)$ 
a homomorphism.
Then $\ro(\Gamma)$ has a finite orbit.
\end{mthm}
%

Theorems~\ref{T:product}, \ref{T:Ghys}, \ref{T:higherrank} 
should be considered as \emph{higher rank superrigidity} with $\Homeo_{+}(\Circle)$
as a target. The following result is an abstract \emph{commensurator superrigidity}.
\begin{mthm}[Commensurator Superrigidity]
\label{T:commensurator}
Let $\Gamma<\Lambda<G$ be groups, where $G$ is a  \lcsc group which has no
infinite discrete quotients, $\Gamma$ is a lattice in $G$,
$\Lambda$ commensurates $\Gamma$ and is dense in $G$.
Let $\ro:\Lambda\to\Homeo_{+}(\Circle)$ be a homomorphism such that
$\ro(\Gamma)$ is non-elementary.

Then $G$ admits an epimorphism $p:G\to\PSL_{2}(\RR)$ so that 
the original action $\ro:\Lambda\to \Homeo_{+}(\Circle)$ is semi-conjugate to
the action 
\[
	\ro_{0}:\Lambda< G\overto{p}\PSL_{2}(\RR)<\Homeo_{+}(\Circle)
\]
by fractional linear transformations. 
If $G$ is topologically simple, then $G\cong \PSL_{2}(\RR)$ and $\Gamma$
is arithmetic.
\end{mthm}

\medskip

There is an interesting class of discrete countable groups -- the so called $\tilde{A}_{2}$ 
groups -- which share many common properties with higher rank lattices, but do not appear
as lattices in higher rank semi-simple groups. 
These groups act simply transitively on thick $\tilde{A}_{2}$ buildings ${\X}$.
If $k$ is a non-Archemidean field with $\mathcal{O}$ denoting the ring of integers in $k$,
$\Gamma=\PGL_{3}(\mathcal{O})$ is such an example -- it acts on the Bruhat-Tits 
building ${\X}$ of $G=\PGL_{3}(k)$.
However, it turns out that there are many other \emph{exotic} 
thick $\tilde{A}_{2}$ buildings which 
admit a simply transitive group of isometries $\Gamma$.
These groups, commonly known as $\tilde{A}_{2}$ \emph{groups}, were constructed and 
described by  Cartwright, Mantero and Steger \cite{CartwrightManteroSteger:A2I:93}, 
\cite{CartwrightManteroSteger:A2II:93}.
We announce the following result, which seems to be the first superrigidity result for such groups.
\begin{mthm} 
\label{T:A2}
Let $\Gamma$ be an $\tilde{A}_{2}$ group as above. 
Then any action $\ro:\Gamma\to\Homeo_{+}(\Circle)$ has a finite orbit.
\end{mthm}
Here we shall outline the proof of this result and its cocycle version 
(Theorem \ref{T:cocycle-A2} below), but will provide the full details in further papers.


\subsection*{Cocycle versions of the results}
We have already mentioned circle bundle version of Ghys' theorem, 
proved by Witte and Zimmer \cite{{WitteZimmer:circlebundle:01}}.
In \cite{Navas:cocycles:06} Navas proves some superrigidity
results for cocycles $\ro:G\times X\to \Diff^{1+\tau}(\Circle)$.
Let us briefly describe the general setting for such results. 

Let $G$ be some \lcsc group, $(X,\mu)$ a standard probability space
with a measurable, measure preserving, ergodic action of $G$.
We consider $G$-actions on \emph{measurable circle bundles} over $X$.
These are trivial bundles and, identifying them with $X\times\Circle$,
a $G$-action corresponds to a measurable cocycle
\[
	\ro:G\times X\to\Homeo_{+}(\Circle),\qquad
	g:(x,p)\mapsto (g.\omega, \ro(g,x).p).
\] 
For $\ro$ to be a cocycle means that for all $g,h\in G$ and $\mu$-a.e. 
$x\in X$
\[
	\ro(g h, x)=\ro(g,h.x)\circ\ro(h,x).
\]
\begin{defns}
\label{D:cocycle}
We shall say that two cocycles $\ro,\ro': G\times X\to \Homeo_{+}(\Circle)$
are \emph{cohomologous}, and that the corresponding $G$-actions are 
\emph{conjugate} if there exists  a measurable map $X\to\Homeo_{+}(\Circle)$, 
$x\mapsto f_{x}$, so that
\[
	\ro'(g,x)=f_{g.x}\circ\ro(g,x)\circ f_{x}^{-1}.
\]
Cocycle $\ro$ will be said to be \emph{semi-conjugate} to $\ro'$, and
the corresponding $G$-actions \emph{semi-conjugate}, if there exists
a measurable family of maps $\{f_{x}:\Circle\to\Circle \mid x\in X\}$
where each $f_{x}$ is a continuous, monotonic, surjection $\Circle\to\Circle$, 
so that 
\[
	\ro'(g,x) \circ f_{x}=f_{g.x}\circ\ro(g,x).
\]
Finally, a $G$-action on such a bundle is \emph{elementary} if it preserves a probability 
measure $\tilde\mu$ on $X\times\Circle$ the projection of which to $X$ is $\mu$.  
Equivalently, a cocycle $\ro$ is elementary if it is cohomologous to a cocycle taking 
values in an elementary subgroup $H_{0}<\Homeo_{+}(\Circle)$.
\end{defns}
The following is a cocycle analogue of Theorem~\ref{T:product} 
(compare with Navas' \cite[Theorem B]{Navas:cocycles:06}):
\begin{mthm}
\label{T:cocycle-product}
Let $G=G_{1}\times\cdots\times G_{n}$ be a product of $n\ge 2$ \lcsc groups,
$(X,\mu)$ a standard probability space with a measurable, measure-preserving
$G$-action in which each factor $G_{i}$ acts ergodically.
Let $\ro:G\times X\to\Homeo_{+}(\Circle)$ be a measurable cocycle,
such that the corresponding $G$-action on $X\times\Circle$ is non-elementary.

Then one of the factors $G_{i}$ admits a continuous epimorphism 
$p:G_{i}\to\PSL_{2}(\RR)$
and the cocycle $\ro$ is semi-conjugate to the homomorphism
\[
	\ro_{0}:G\overto{\pi_{i}}G_{i}\overto{p}\PSL_{2}(\RR)<\Homeo_{+}(\Circle).
\]
\end{mthm}

Let us state a result of Witte and Zimmer, which is 
a cocycle generalizations of Ghys' theorem (though only of the 
case where $G$ has no $\PSL_{2}(\RR)$ factors). 
We state the case of a simple $G$:
\begin{thm}[Witte and Zimmer \cite{WitteZimmer:circlebundle:01}]
\label{T:WZ}
Let $G=\mathbf{G}(k)$ be a higher rank group as in Theorem~\ref{T:higherrank}.
Then for any ergodic action of $G$, or of a lattice $\Gamma<G$, on a probability measure
$(X,\mu)$ all measurable cocycles $\ro:G\times X\to\Homeo_{+}(\Circle)$, 
or $\ro:\Gamma\times X\to\Homeo_{+}(\Circle)$, are elementary.
\end{thm}
\begin{rem}
Note that Theorem~\ref{T:cocycle-product}, in particular, covers
the case where higher rank semi-simple $G$ has a $\PSL_{2}(\RR)$ factor/s.
This completes the cocycle generalization of Ghys' result, 
left open in \cite{WitteZimmer:circlebundle:01}. 
\end{rem}

Finally we state the cocycle version of Theorem~\ref{T:A2}.
\begin{mthm}
\label{T:cocycle-A2}
Let $\Gamma$ be an $\tilde{A}_{2}$ group, acting ergodically on a probability 
space $(X,\mu)$.
Then all measurable cocycles $\ro:\Gamma\times X\to\Homeo_{+}(\Circle)$
are elementary.
\end{mthm}

\medskip

\subsection{Remarks on Elementary actions} \label{SS:elementary}

Recall that most groups $\Gamma$ discussed above
have trivial $H^{1}(\Gamma,\RR)$, i.e. have no homomorphisms to $\RR$.
Specifically this is known for the following classes of groups:
\begin{enumerate}
\item
	irreducible lattices in higher rank semi-simple algebraic groups $G$
	(c.f. Margulis \cite{Margulis:book:91});
\item
	irreducible cocompact lattices $\Gamma$ in products 
	$G=G_{1}\times\cdots\times G_{n}$ of $n\ge 2$ general \lcsc groups,
	with $H^{1}(G_{i},\RR)=\{0\}$ for all $1\le i\le n$ 
	(Shalom \cite[Theorem 0.8]{Shalom:Inven:00});
\item
	dense commensurators $\Lambda$ of a cocompact lattice  $\Gamma<G$ in a 
	general \lcsc groups $G$ with $H^{1}(G,\RR)=\{0\}$
	(Shalom \cite[Theorem 0.8']{Shalom:Inven:00}).
\item
	$\tilde{A}_{2}$ groups, because they have Kazhdan's property (T) 	
	(\cite{CartwrightMlotkowskiSteger:T:94}).
\end{enumerate}
Hence elementary actions of such groups are those which have a finite orbit.
The restriction of such an action to an appropriate subgroup of finite index
will have a fixed point. These subgroups still have vanishing $H^{1}$.

It was pointed out by Ghys \cite{Ghys:lattices:99} 
(and by Burger and Monod \cite{BurgerMonod:JEMS:99}) that if such a group acts by 
$C^{1}$-diffeomorphisms, i.e. $\ro(\Gamma)<\Diff^{1}_{+}(\Circle)<\Homeo_{+}(\Circle)$ then, 
it follows from Thurston's stability theorem \cite{Thurston:Reeb:74}
that the action of the subgroup is trivial (and the action of the original one is finite).
Without the $C^{1}$-smoothness assumption it is, in general, open
whether a group as above can have an infinite action with a finite orbit on the circle.
Equivalently, it is not known whether some finite index subgroup of a group on the above list
can have an infinite action on an interval: $\Gamma\to\Homeo_{+}([0,1])$. 
It is widely believed that such actions do not exist for higher rank lattices.
This has been proved by Dave Witte \cite{Witte:circle:94} for the case
$\text{rank}_{\QQ}(\Gamma)\ge 2$, and more recently by 
Lucy Lifschitz and Dave (Witte) Morris \cite{MorrisLifschitz:Qrank1:04}
for many $\text{rank}_{\QQ}(\Gamma)=1$ cases. 
The problem is still open for all cocompact (i.e. $\QQ$-rank $0$)
higher rank lattices.

\subsection{Remarks on Semi-conjugacy}
The semi-conjugacy $f$ in definition \ref{D:semi-conj} can always 
be presented as a composition $f=g\circ h$ of a degree one, non-decreasing, 
continuous surjection $h:\Circle\to\Circle$, and a degree $k\ge 1$ covering map 
$g(z)=z^{k}$.
If the original action is assumed to be minimal, then $h$ is a homeomorphism.
But even in this case semi-conjugacy cannot be replaced by a conjugacy 
in the statements of Theorems~\ref{T:product}, \ref{C:arithm}, \ref{T:commensurator}. 
The basic counterexample is an irreducible lattice $\Gamma$ in
$G_{0}=\PSL_{2}(\RR)\times\PSL_{2}(\RR)$ which has a lifted embedding into
\[
	G=\PSL_{2}(\RR)\times\PSL_{2}(\RR)_{k}
\] 
where $\PSL_{2}(\RR)_{k}$ is an $\ZZ/k\ZZ$ central extension of $\PSL_{2}(\RR)$
corresponding to the $k$-fold cover $\Circle\to\Circle$.
The $\Gamma$-action $\ro$ via the second factor of $G$ is not
conjugate, albeit semi-conjugate, to the action $\ro_{0}$ via $G$.
Similar examples exist for commensurators.
However assuming that the original action is minimal and strongly proximal 
(see Definition \ref{D:str-prox}) the semi-conjugacies in the above stated results  
can be replaced by conjugacies.

\subsection*{Organization of the paper}
The rest of the paper is organized as follows.
In section~\ref{S:Boundary} we discuss all the concepts and results which reflect
the ``higher rank'' phenomena of the acting groups.
This includes the discussion of the \emph{boundaries} that we use, the concept
of the generalized Weyl group, basic examples and basic constructions, such as the
Galois correspondence between subgroups of the Weyl and quotients of the boundary.
Section~\ref{S:Circle} is focused on the target group $\Homeo_{+}(\Circle)$:
we recall some properties of group dynamics on the circle
and establish some special properties related to the boundary theory 
(Theorem~\ref{T:reduction2atoms}) and the generalized Weyl group (Theorem~\ref{T:cyclic}).
These results are put together in Section \ref{S:mainproofs} to prove
Theorems \ref{T:product}, \ref{T:higherrank}, \ref{T:A2}. 
Theorem \ref{T:commensurator}, which is somewhat independent, is proved in 
Section \ref{S:commensurator}.
Section \ref{S:cocycles} contains the cocycle versions of the results,
with some of the technical discussion postponed to Appendix \ref{S:Appendix}.

\subsection*{Acknowledgments}
We would like to thank Dave (Witte) Morris for his inspiring lectures about lattice actions on
the circle given in Chicago in the spring 2005. 
We also thank Andr\'es Navas for his helpful discussions about
\ref{T:Margulis-Ghys}.

\subsection*{Dedication} 
Ali Shaker has tragically died in September 2005. 
At the time Ali was working on questions of commensurator superrigidity, and
Theorem~\ref{T:commensurator} was supposed to be part of his PhD thesis.
This paper is dedicated to his memory.

\section{Boundary Theory and Weyl Groups}\label{S:Boundary}

\subsection{$G$-boundaries}

Boundary Theory is a broad and somewhat vague term which describes ways
to associate a $G$-space -- a $G$-\emph{boundary} -- to a given group $G$ in a way
which helps to understand the structure of the group $G$ itself and of its actions
on other spaces.
In this paper we shall use the concept of \emph{amenable, doubly ergodic with unitary 
coefficients} boundaries, introduced by Burger and Monod in the course of
their study of bounded cohomology (we shall use a slightly weaker variant, see below).
The main point is that this notion of a boundary is flexible enough to pass to lattices
(see Proposition~\ref{prop:boundary-lattice}) 
and to be applicable to measurable $G$-cocycles, 
and yet it is sufficiently powerful for our applications. 
Also every \lcsc group $G$ admits a $G$-boundary in this sense, even though not a unique one.
\begin{defn}
\label{D:G-boundary}
Let $G$ be an arbitrary \lcsc group. 
Hereafter a $G$-\emph{boundary} means a Lebesgue space $(B,\nu)$ with a measure class preserving measurable $G$-action $G\times B\to B$ which is
\begin{description}
\item[Amenable] in the following sense: for any continuous 
	$G$-action $G\to\Aff(Q)$ by affine maps on a compact convex
	subset $Q$ in some locally convex topological vector space,
	there exists a measurable $G$-equivariant map $\Phi:B\to Q$.
\item[Doubly Ergodic] in the following sense:
	for every Lebesgue space, $(X,\xi)$, with a {\em measure preserving} $G$ action
	$G \times X \to X$, the projection map $p:(B\times B \times X,\nu \times \nu\times \xi) \to (X,\xi)$
	gives rise to a Lebesgue-space isomorphism at the level of ergodic componenets:
	$(B \times B \times X) \mod G \overto{\sim} X \mod G$\footnote{Given a measurable, 
	measure class preserving action of a \lcsc group $G$ on a standard Lebesgue measure space 
	$(X,\xi)$ we denote by $X\mod G$ the Lebesgue space of ergodic components of the action.}.
\end{description}
\end{defn}
%
The amenability condition is just slightly weaker than Zimmer's original definition 
\cite{Zimmer:Amenable:78} (which allows $Q$ to vary measurably over $x\in B$).
The double ergodicity condition is equivalent to the condition "For every
measure preserving ergodic $G$-space $X$, $B \times B \times X$ is ergodic" by a use of 
an apropriate ergodic decomposition.
This condition follows from the original condition of {\em double ergodicity with unitary coefficients} 
as in Burger-Monod \cite[Definition 5, p. 221]{BurgerMonod:GAFA:02}:
%
\begin{lem} \label{unitary-coeff}
Assume the Lebesgue $G$-space $X$ is such that for every unitary 
$G$-representation $\pi$ on $\sH$, the obvious map (given by constant functions) 
$\sH^{\pi(G)} \rightarrow L^{\infty}(X,\sH)^G$ is an isomorphism
($X$ is {\em ergodic with unitary coefficients}). 
Then for every probability space, $(Y,\eta)$, with a {\em measure preserving} 
$G$-action, the map
$(X \times Y) \mod G {\to} Y \mod G$
is an isomorphism.
\end{lem}

\begin{proof}
Consider the following, naturally given, diagram
\[  \begin{array}{ccccc}
    L^{\infty}(Y)^G & \to & L^{\infty}(X,L^{\infty}(Y))^G &
    \overset {\sim}{\to} & L^{\infty}(X\times Y)^G \\
    \cap & & \cap & & \\
    L^2(Y)^G & \to & L^{\infty}(X,L^2(Y))^G & &
    \end{array}
\]
Note that if the map $L^2(Y)^G \to L^{\infty}(X,L^2(Y))^G$ is an isomorphism then so is
the composed map $L^{\infty}(Y)^G \to L^{\infty}(X\times Y)^G$.
The latter gives a Lebesgue isomorphism between the corresponding von-Neumann spectra.
\end{proof}
We have already mentioned that every \lcsc group $G$ admits a boundary $(B,\nu)$ 
in the above sense.
Haven singled out the key properties of amenability + double ergodicity with coefficients,
Burger and Monod have proved in \cite[Theorem 6]{BurgerMonod:GAFA:02}
that every compactly generated \lcsc group $G$ admits a boundary in the above sense.
Shortly afterwards Kaimanovich removed the compact generation assumption,
by an elegant use of random walks, namely using Furstenberg's notion of the 
Poisson boundary of a random walk on a group.
\begin{thm}[Kaimanovich \cite{Kaimanovich:DE:03}]
\label{T:Kaim}
Let $\mu$ be an \'etal\'ee symmetric probability measure on a \lcsc group $G$,
and let $(B,\nu)$ be the associated Poisson boundary of $(G,\mu)$.
Then $(B,\nu)$ is amenable and double ergodic with unitary coefficients,
and in particular it is a $G$-boundary in the above sense.
\end{thm}

\medskip

Let us record some useful observations. The following is immediate from the definitions:
\begin{prop}[Amenability Criterion]
\label{P:amenability}
For a non-amenable \lcsc group $G$ all its boundaries $(B,\nu)$
have the measure type of the interval, i.e. they are Lebesgue non-atomic measure spaces.
For an amenable $G$ the trivial space -- the singleton -- may serve as a boundary.
\end{prop}
\begin{rem}
Note that for an amenable $G$ any weakly mixing \emph{measure preserving}
action gives a $G$-boundary. Such actions can have huge $W_{G,B}$ (which will
never produce any rigidity phenomena in our sense). 
Some amenable groups might have a ``genuine'' boundary with 
a non-trivial Weyl group, e.g. Kaimanovich and Vershik showed that any 
symmetric random walk with first moment on  the wreath product 
$G=\ZZ^{3}\ltimes\sum_{\ZZ^{3}}\ZZ/2\ZZ$ has a non-trivial Poisson boundary
\cite[Example 6.2]{KaimanovichVershik:AnnProb:83}.
\end{rem}
\begin{prop}[Lattices]
\label{prop:boundary-lattice}
Let $\Gamma<G$ be a lattice,
and $(B,\nu)$ be a $G$-boundary.
Then $(B,\nu)$ is a $\Gamma$-boundary.
\end{prop}

\begin{proof}
The amenability of $G$ on $B$ is inherited by any closed subgroup (\cite[Proposition 4.1.6]{Zimmer:book:84}).
Let $(X,\mu)$ be a probability space with an ergodic measure preserving action of $\Gamma$.
Then $G$ acts on the probability space $(G\times X)/\Gamma\cong G/\Gamma\times X$ ergodically.
Hence the $G$-action on $B\times B\times (G\times X)/\Gamma$ is ergodic, which is equivalent
to the ergodicity of the diagonal $\Gamma$-action on $B\times B\times X$. 
\end{proof}

\begin{prop}[Products]
\label{P:products-of-boundaries}
Let $G_1,\ldots,G_n$ be \lcsc groups, and let 
$(B_i,\nu_{i})$ be $G_i$ boundaries.
Set $G=G_{1}\times\cdots\times G_{n}$ and 
$(B,\nu)=(B_1\times\cdots B_n,\nu_{1}\times\cdots\nu_{n})$
with the $G$-action on $B$ given by $(g.b)_{i}=g_{i}.b_{i}$. 
Then $B$ is a $G$-boundary.
\end{prop}

\begin{proof}
Assume $n=2$ (the general case will follow by induction).

{\em Amenability:}
Let $V$ be a locally convex $G$-topological vector space, 
and let $Q \subset V$ be a compact convex $G$-invariant set.
Set $V_1=L^{\infty}(B_1,V)^{G_1}$ and $Q_1=\{f \in V_1~|~f(B_1)\subset Q\}$.
By the $G_1$-amenability of $B_1$, $Q_1$ is not empty.
$G_2$ acts on $V_1$ via its action on $V$, and $Q_1$ is $G_2$-invariant, as $Q$ is.
By the $G_2$-amenability of $B_2$ there is a $G_2$ map from $B_2$ to $Q_1$.
The latter can be interpreted as a $G$-map from $B$ to $Q$.

{\em Double ergodicity} 
Denote $\Omega_i=B_i \times B_i$, and $\Omega=\Omega_1\times \Omega_2$.
Let $X$ be a Lebesgue $G$-space with an invariant measure.
Then
\begin{eqnarray*}
	(\Omega \times X)\mod G&=&(\Omega_2 \times ((\Omega_1 \times X) \mod G_1)) \mod G_2\\
	&=&(\Omega_2 \times (X\mod G_1))\mod G_2 = (X\mod G_1)\mod G_2=X \mod G.
\end{eqnarray*}
\end{proof}

A well known class of examples of $G$-boundaries is given by the maximal 
Furstenberg boundaries of semi-simple Lie-groups. 
A standard example of this sort for $G=\SL_{n}(\RR)$, 
is the space of \emph{flags} -- sequences of nested linear subspaces 
\[
	E_{1}\subset E_{2}\subset\dots\subset E_{n-1}\qquad \text{with}
	\qquad\dim E_{i}=i.
\]
It is the compact homogeneous $G$-spaces $G/P$ where $P$ is the upper triangular subgroup.
More precisely, and more generally, we have the following.
\begin{lem}[Boundaries of Semi-simple Groups]
\label{L:semi-simple}
Let $k$ be a local field and $\mathbf{G}$ be a semi-simple algebraic group defined over $k$.
Denote by $\mathbf{P}$ a minimal parabolic subgroup of $\mathbf{G}$
defined over $k$, and by $\mathbf{S}$ a maximal $k$-split torus in $\mathbf{P}$.
Endowed with its Hausdorff topology $G=\mathbf{G}(k)$
is a locally compact group, $P=\mathbf{P}(k)$ and $S=\mathbf{S}(k)$
are closed subgroups of $G$. The homogeneous $G$-space $B=G/P$
with the Haar measure class $\nu$ is a $G$-boundary.
\end{lem}
Before proving the lemma, let us understand better the space $B\times B$.

\begin{lem} \label{L:BxB=G/Z}
The spaces $(B\times B,\nu \times \nu)$ and $G/Z_G(S)$ are
isomorphic as $G$-Lebesgue spaces, i.e there is an invertable
measur class preserving $G$-equivariant map defined up to measure zero
between them (where $G/Z_G(S)$ is taken with the Haar measure class).
\end{lem}

\begin{proof} 
By the Bruhat decomposition there are finitly many $G$-orbits in $B\times B$
(parametrized by the $k$-relative Weyl group).
All of those orbits but the unique open one, 
are lower dimensional in the $k$-analityc sense, 
and therefore have measure zero.
Thus, $B\times B$ minus a null set equals the unique open orbit in it,
namely the orbit in which the stabilizers are reductive.
By \cite[4.8]{BorelTits:groupsreductifs:65}, the point stabilizers of this orbit are conjugates of
$Z_{\mathbf{G}}(\mathbf{S})(k)=Z_G(S)$.
Thus, $B\times B$ is (up to a null set) isomorphic to $G/Z_G(S)$ as a $G$-space.
\end{proof}

\begin{proof}[Proof of Lemma~\ref{L:semi-simple}]
The amenability of the homogeneous $G$ space $G/P$ is equivalent to the
amenability of the stabilizer $P$ of a point.
We have $P=\mathbf{P}(k)=Z_{\mathbf{G}}(\mathbf{S})(k)\ltimes R_u(\mathbf{P})(k)$.
The unipotent radical $R_u(\mathbf{P})(k)$ is amenable being nilpotent,
and $Z_{\mathbf{G}}(\mathbf{S})(k)$ is amenable being a compact exstension of an abelian group
(the fact that $Z_{\mathbf{G}}(\mathbf{S})(k)/\mathbf{S}(k)$ is compact is in \cite[Proposition (2.3.6)]{Margulis:book:91}). Hence the $G$-action on $G/P$ is amenable.

Double ergodicity with unitary coefficients for $G/P$ follows from the Howe-Moore theorem.
Any measurable $G$-equivariant map from $B\times B \simeq G/Z_G(S)$ to a Hilbert space sends
the coset $Z_G(S)$ to a $Z_G(S)$ fixed point $v$.
By Howe-Moore $v$ is $\mathbf{G}^+(k)$ invariant.
Then $v$ is $G$-invariant, because
$G=\mathbf{G}(k)^{+}\cdot Z_\mathbf{G}(\mathbf{S})(k)$
\cite[Proposition 6.11(i)]{BorelTits:hom-abs:73}.
We are done by Lemma~\ref{unitary-coeff}.
\end{proof}
\begin{rem}
Another way of proving lemma \ref{L:semi-simple}
is by using the fact that $G/P$ (with the Haar measure class) can be realized as 
the Poisson boundary of $G$
(with respect to any admissible symmetric measure).
See \cite{Furstenberg:Poisson:63} for real Lie-groups,
\cite{GuivarchJiTaylor:book:98} for other local fields, simply connected groups,
and \cite[Proposition 5.2]{BaderShalom:Factor:06} for the general case.
\end{rem}

\medskip

Let ${\X}$ be an $\tilde{A}_{2}$ building. 
Let us briefly describe the notion of a geometric boundary $B$ of ${\X}$.
${\X}$ is a union of \emph{apartments}, which are planes divided by three 
families of parallel lines, the \emph{walls}, into equilateral triangles,  the \emph{chambers}.
In a given apartment any vertex is contained in $3$ walls which divide it into $6$ connected
components, called sectors. 
Two sectors are declared to be equivalent if there exists a sector contained in both.
$B$ is the space of equivalence classes of sectors in ${\X}$.
It carries a natural topology which makes it into a compact (Cantor) space.
The group of isometries of ${\X}$ acts on $B$ and this action is continuous.
We consider the special case where the $\tilde{A}_{2}$ building ${\X}$
has a group $\Gamma$ of isometries which acts transitively on vertices. 
In \cite{ManteroZappa:00} (see also \cite{Cartwright:Harmonic:99})
Mantero and Zappa showed a one-to-one correspondence between 
harmonic functions on $\X$ and finitely additive measures on $B$ (given by the Poisson transform). 
One can deduce from this result that the Poisson boundary of an isotropic random walk on $\Gamma$
can be realized 
as some probability measure on $B$. 
Hence applying Kaimanovich theorem \ref{T:Kaim} we obtain:
%
\begin{prop}[Boundary for $\tilde{A}_{2}$ groups]
\label{P:A2bnd}
There is a probability measure $\nu$ on the geometric boundary $B$ of ${\X}$
so that $(B,\nu)$ is a $\Gamma$-boundary.
\end{prop}


\subsection{Generalized Weyl Groups}
\label{SS:Weyl}
\begin{defn}
\label{D:Weyl}
Let $G$ be a locally compact group, $(B,\nu)$ a $G$-boundary.
The associated \emph{Weyl group} $W_{G,B}$ is defined to be 
the group of all measure class preserving automorphisms
of $B\times B$ commuting with the diagonal $G$-action.

If $B$ is non-trivial, the involution $w_{0}:(x,y)\mapsto (y,x)$ of $B\times B$
is a non-trivial element of $W_{G,B}$. 
We shall refer to it as \emph{the long element} of the Weyl group.
\end{defn}

For a semi-simple group $G$ and its Furstenberg boundary $B=G/P$
the Weyl group $W_{G,B}$ coincides with the classical Weyl group $W$
of $G$, and $w_{0}$ corresponds to \emph{the long element} of $W$.
More precisely, we have the following.
\begin{prop} \label{P:W4ss}
Let $G=\mathbf{G}(k)$ be the $k$-points of a semi-simple group defined over $k$
and $B=G/P$ be as in Lemma~\ref{L:semi-simple}. 
Then the generalized Weyl group $W_{G,B}$ coinsides with the $k$-relative Weyl group, 
$W=N_{\mathbf{G}}(\mathbf{S})/Z_{\mathbf{G}}(\mathbf{S})$.
\end{prop}
Given a group $G$, a subgroup $H$, let $\Aut_{G}(G/H)$ denote the group of 
permutations of the set $G/H$ commuting with the $G$-action. 
If $G$ is a \lcsc group and $H$ is closed, $G/H$ carries a unique 
$G$-invariant measure class $[m_{G/H}]$.
We denote by $\Aut_{G}(G/H,[m_{G/H}])$ 
the group of all measure class preserving measurable bijections of $G/H$,
modulo the equivalence relation of being a.e. the same.
We need the following simple. but useful 
\begin{lem} \label{L:Aut=N/H}
There is a natural isomorphism $\Aut_{G}(G/H)\cong N_{G}(H)/H$; if 
$G$ is a \lcsc group and $H<G$ is closed, then these groups are also 
naturally isomorphic to $\Aut_{G}(G/H,[m_{G/H}])$.
\end{lem}
\begin{proof}
First observe that the map $\alpha:N_G(H)\to\Aut_{G}(G/H)$, 
$\alpha(n): gH\mapsto gnH$, 
is well defined, defines a group homomorphism, and $\Ker(\alpha)=H$.
Hence it gives an embedding $\bar{\alpha}:N_{G}(H)/H\to\Aut_{G}(G/H)$.
To see that it is a surjection, note that given $\phi\in\Aut_{G}(G/H)$ 
any $g_{0}\in G$ with $\phi(eH)=g_{0}H$ satisfies
\[
	g_{0}Hg_{0}^{-1}=\Stab(g_{0}H)=\Stab(\phi(eH))=\Stab(eH)=H. 
\]
Thus $g_{0}\in N_G(H)$. Moreover $\phi(gH)=\phi(geH)=g\phi(eH)=gg_{0}H$
for any $g\in G$, that is $\phi=\alpha(g_{0})$.

In the context of \lcsc groups, elements of $N_{G}(H)/H$ still define 
measurable centralizers of the $G$-action on $(G/H,[m_{G/H}])$,
and it is easy to see that non-trivial elements of $N_{G}(H)/H$ remain non-trival 
even measurable. 
On the other hand, if $\phi:G/H\to G/H$ is measurable
and for each $g\in G$, $\phi(g.x)=g.\phi(x)$ for a.e. $x\in G/H$; then
by Fubini for some (in fact a.e.) $x=g_{0}H\in G/H$ and for a.e. $g\in G$ we have
$g^{-1}\phi(gg_{0}H)=\phi(g_{0})$. Denoting $\phi(g_{0}H)=g_{1}H$ and 
$n=g_{0}^{-1}g_{1}$ we get $\phi(g'H)=g' n H$ for all $g'$ in some conull set  $E\subset G$.
For arbitrary $h\in H$ the set $E\cap Eh^{-1}$ is conull in $G$, and for any $g\in E\cap Eh^{-1}$: 
\[
	gn H=\phi(g H)=\phi(gh H)=gh n H
\]
giving $n^{-1}hn\in H$. Thus $n\in N_{G}(H)$ and $\phi$ a.e. agrees with $\alpha(n)$.
\end{proof}

\begin{proof}[{Proof of Proposition~\ref{P:W4ss}}]
By Lemmas~\ref{L:BxB=G/Z} and~\ref{L:Aut=N/H},
the abstract Weyl group $W_{G,B}$ can be identified with $N_G(Z_G)/Z_G$.
As both $\mathbf{S}$ and $Z_{\mathbf{G}}(\mathbf{S})$ are reductive,
each is equal to the Zariski closure of its $k$-points
\cite[2.14(c)]{BorelTits:groupsreductifs:65}.
It follows that $Z_G(S)=Z_{\mathbf{G}}(\mathbf{S})(k)$,
and $N_G(Z_G(S))=N_{\mathbf{G}} (Z_{\mathbf{G}}(\mathbf{S}))(k)$.
Since $\mathbf{S}$ is the unique maximal $k$-split torus in
$Z_{\mathbf{G}}(\mathbf{S})$, we get also
$N_{\mathbf{G}} (Z_{\mathbf{G}}(\mathbf{S}))=
N_{\mathbf{G}} (\mathbf{S})$.
This gives,
$N_G(Z_G(S))=N_{\mathbf{G}} (\mathbf{S})(k)$ and, 
using \cite[5.4]{BorelTits:groupsreductifs:65}, allows to conclude
\[
	W_{G,B} \cong N_G(Z_G(S))/Z_G(S) \cong 
	N_{\mathbf{G}} (\mathbf{S})(k)/Z_{\mathbf{G}}(\mathbf{S})(k)
	\cong N_{\mathbf{G}}(\mathbf{S})/Z_{\mathbf{G}}(\mathbf{S})=W.
\]
\end{proof}

\begin{obs}
\label{O:Weyl4products}
Let $G=G_{1}\times\cdots\times G_{n}$ be a direct product of $n$ general 
\emph{non-amenable} locally compact groups, and 
$
	(B,\nu)=\prod_{i=1}^{n} (B_{i},\nu_{i})
$
as in Proposition~\ref{P:products-of-boundaries}.
Then $W_{G,B}$ contains $\prod W_{G_{i},B_{i}}$ with each 
$W_{G_{i},B_{i}}$ being non-trivial. 
In particular, for every $i=1,\ldots,n$ we have an involution in $W_{G,B}$ given by
\[
	\tau_{i}: (x_{1},\dots,x_{i},\dots,y_{1},\dots,y_{i},\dots)\ \mapsto
	(x_{1},\dots,y_{i},\dots,y_{1},\dots,x_{i},\dots)
\]
where $x_{j},y_{j}\in B_{j}$.
The group generated by these involutions $\tau_i$ is isomorphic to
$(\ZZ/2\ZZ)^{n}$. We denote it by $W$.
For a subset $I$ of $\{1,\ldots,n\}$, we denote $\tau_I=\prod_{i\in I}\tau_i$,
and $W_I=\prod_{i\in I} \langle \tau_i \rangle=\{ \tau_J~|~J\subset I\}$.
The \emph{long element} of $W_{G,B}$ is given by $w_{0}=\tau_{\{1,\ldots,n\}}$.
\end{obs}

\begin{prop}
\label{P:W4A2}
Let $\Gamma$ be an $\tilde{A}_{2}$ group, and $(B,\nu)$ the boundary
as in \ref{P:A2bnd}. Then $W_{G,B}$ contains $S_{3}$.
\end{prop}
\begin{proof}[Sketch of the Proof]
For $\nu\times\nu$-a.e. two equivalence classes $\sigma, \sigma'\in B$ of sectors 
there exists a unique apartment $P$ which contains representatives $S\in\sigma$, $S'\in\sigma'$.
which form opposite sectors in $P$, and therefore, define a pair of opposite
sides in the \emph{hexagon at infinity} $\pa P$ of $P$. 
The group $S_{3}$ acts on such objects by preserving each of the planes $P$ but permuting
the choice of the three pairs of sides on the hexagon $\pa P$ (it is generated by the reflections 
around the ``main diagonals'' of the hexagon).
Being geometrically defined this $S_{3}$-action commutes with the group of isometries $\Gamma$;
it also preserves the measure class of $\nu\times\nu$,
\end{proof}

\subsection{The Galois Correspondence}
\label{S:Galois}
The abstract term \emph{Galois correspondence} (often called Galois \emph{connection})
refers to a situation where two 
posets (partially ordered sets) $\mathcal{A}$ and $\mathcal{B}$ are related
by two order reversing maps $\alpha:\mathcal{B}\overto{}\mathcal{A}$ and
$\beta:\mathcal{A}\overto{}\mathcal{B}$, so that
\[
	\beta(A)\le B\ \Longleftrightarrow\  A\le \alpha(B),
	\qquad\qquad(A\in\mathcal{A},\ B\in\mathcal{B}).
\]
In this general setting, it can be shown that 
$\alpha\circ \beta:\mathcal{A}\to\mathcal{A}$ and 
$\beta\circ\alpha:\mathcal{B}\to\mathcal{B}$ are order preserving idempotents.
This enables one to define the corresponding \emph{closure} operations:
\[
	\ol{A}=\alpha\circ\beta(A),\qquad \ol{B}=\beta\circ\alpha(B).
\]
The collections $\ol{\mathcal{A}}$ and $\ol{\mathcal{B}}$ of closed objects 
in $\mathcal{A}$ and $\mathcal{B}$ (which can be viewed both as sub posets or as 
quotient posets of $\mathcal{A}$ and $\mathcal{B}$) are order anti-isomorphic via
the restrictions/quotients $\bar\alpha=\bar\beta^{-1}$ of $\alpha$ and $\beta$.
The classical Galois correspondence relates subfields to subgroups.
We shall consider the following instance of the abstract Galois correspondence.

Let $G$ be a \lcsc group, $(B,\nu)$ be a $G$-boundary, and $W<W_{G,B}$ 
be some closed subgroup of the corresponding abstract Weyl group. 
(In what follows $W$ will always be finite and will be assumed to contain the 
long element $w_{0}$).
Consider the poset $\mathcal{SG}(W)$ consisting of all subgroups of $W$ 
ordered by inclusion.

Denote by $\mathcal{Q}(B)$ the poset of all measurable quotients of $(B,\nu)$. 
Elements of $\mathcal{Q}(B)$ are measurable 
quotient maps $\Phi:(B,\nu)\overto{}(C,\eta)$ with $\Phi '\le \Phi$
(or $C'\le C$) if $\Phi':B\overto{} C'$ factors through $\Phi:B\overto{} C$.
Alternatively, one might  think of elements of $\mathcal{Q}(B)$ as 
complete sub $\sigma$-algebras $\mathcal{L}_{C}$ of the $\sigma$-algebra 
$\mathcal{L}_{B}$ of $B$, or as 
von Neumann subalgebras of $L^{\infty}(B,\nu)$,
ordered by inclusion.

We think of $\mathcal{Q}(B)$ as imbedded in $\mathcal{Q}(B\times B)$,
say via $p_{1}^{*}:L^{\infty}(B)\to L^{\infty}(B\times B)$, $(p_{1}^{*}f)(x,y)=f(x)$,
corresponding to the projection $p_{1}:B\times B\to B$ to the first factor.
$W$ has a natural action on $L^{\infty}(B\times B)$, and on $\mathcal{Q}(B\times B)$.
Consider the maps
\[ 
	\mathcal{SG}(W) \overto{T}\mathcal{Q}(B),\qquad 
	\mathcal{SG}(W) \overfrom{S} \mathcal{Q}(B)
\]
defined as follows.
Given a quotient $\Phi:B\to C$ the associated subgroup of $W$ is:
\[
	S(\Phi)=\left\{w\in W\ \left|\ \Phi \circ p_{1}\circ w=\Phi\circ p_{1}\right.\right\}.
\]
Equivalently,
$S(\Phi)=\left\{w\in W\ \left|\ 
	w(E\times B)\simeq  E\times B,\quad \forall E\in \mathcal{L}_{C}\right.\right\}
$,
where $\simeq$ denotes equality modulo $\nu\times\nu$-null sets.
The map $T:\mathcal{SG}(W)\to\mathcal{Q}(B)$ is given by
\[
	T(W')=\left\{E\in\mathcal{L}_{B}\ \left|\ w(E\times B)\simeq E\times B,\ \forall w\in W'
	\right.\right\}.
\]
Note that $T(W')$ is always a complete sub $\sigma$-algebra of $\mathcal{L}_{B}$.
Since $W$ centralizers $G$ this $\sigma$-algebra is $G$-invariant, and therefore 
defines a $G$-equivariant quotient of $B$.
Thus for any $\Phi:B\to C$ in $\mathcal{Q}(B)$ the \emph{closure} 
$\ol{\Phi}=T\circ S (\Phi)$ is a $G$-equivariant quotient $\bar{\Phi}:B\to\bar{C}$
which is finer ($\ge$) than $\Phi$, i.e., for some $F:\bar{C}\to C$ we have
\[
 	\Phi:B\overto{\bar{\Phi}}\bar{C}\overto{F} C.
\]

For a general quotient $\Phi:B\to C$ one expects  $W'=S(\Phi)$ to be the trivial
subgroup $\{\Id_{B\times B}\}$,
in which case the closure of $\Phi:B\to C$ would be the identity map 
$\bar{\Phi}:B\to \bar{B}=B$. 
In this situation the Galois correspondence does not provide any new information.

However, if the subgroup $W'=S(\Phi)$ corresponding to a quotient $\Phi:B\to C$ 
is non-trivial, or more precisely if the $W$-action on $W/S(\Phi)$ has a kernel, 
then the closure $\bar{\Phi}:B\to\bar{C}$ becomes an interesting object,
as examples below show.
Let us first point out a construction which will provide an information about $S(\Phi)$ 
in our setting.
\begin{cons} \label{cons:vec-map}
Let $\Gamma<G$ be a lattice, and $\Phi:B\overto{} C$ be a $\Gamma$-equivariant 
quotient map.
To each element $w\in W/S(\Phi)$ we associate the map 
\[
	\Psi^{(w)}:B \times B \overto{w} B\times B\overto{p_{1}} B\overto{\Phi} C,
\]
and construct the joining 
\[
	\vec{\Psi}:B\times B \overto{} \prod_{W/S(\Phi)} C.
\]
This map is $\Gamma\times W$-equivariant, where 
$\Gamma$ acts diagonally on the target, while $W$ permutes its factors.
\end{cons}
We shall see in Theorem~\ref{T:cyclic} that for boundary maps $\Phi:B\to\Prob(\Circle)$
associated with actions on the circle the $W$-action on $W/S(\Phi)$ necessarily 
factors through a finite cyclic group.

\medskip

Let us now illustrate what implications such information can have on $\Phi$ 
via its closure $\bar{\Phi}$.
\begin{lem} \label{l:long-element}
For an arbitrary \lcsc group $G$, its boundary $B$ and $W<W_{G,B}$ 
containing the long element $w_{0}$, 
the corresponding quotient $T(\langle w_0 \rangle)\in\mathcal{Q}(B)$ is the trivial 
(one point) quotient, and the closure $\overline{\langle w_0 \rangle}$ is all of $W$.
\end{lem}
\begin{proof}
Follows from the definitions.
\end{proof}
\begin{lem} 
\label{L:WPhi4products}
Let $G$ be a product $G=\prod_{i=1}^{n}G_{i}$ of $n\ge 2$ general non-amenable \lcsc groups, 
and that $(B,\nu)$ be a product boundary $\prod_{i=1}^{n}(B_{i},\nu_{i})$ 
as in \ref{O:Weyl4products}.
Assume $W\cong(\ZZ/2\ZZ)^{n}$ is the subgroup generated by the involutions $\tau_i$.
Then
\[
	T\left(\langle \tau_I \rangle\right)=\prod_{j\notin I} B_j, 
	\qquad S(\prod_{j\notin I} B_j)=W_I
	\qquad\text{and so}\qquad
	\overline{\langle \tau_I \rangle}=W_I.
\]
Let $\Phi:B\to C$ be a measurable quotient and suppose that $W_{I}$ is
in the kernel of the $W$-action on $W/S(\Phi)$.
The $\Phi$ factors as follows
\[
	\Phi:B\overto{}\prod_{j\notin I} B_{j}\overto{} C.
\]
\end{lem}
\begin{proof}
Follows from the definitions.
\end{proof}

\begin{prop} 
\label{P:key4simple}
Let $G=\mathbf{G}(k)$ be the $k$-points of a simple algebraic group defined 
over a local field $k$ and $B=G/P$ be as in Lemma~\ref{L:semi-simple}.
Then for any quotient $\Phi:B\to C$ the $W$-action on $W/S(\Phi)$ is faithful,
unless $\Phi$ is the constant map to the singleton.
\end{prop}
The proof of this key Proposition consists of the following two Lemmas:

\begin{lem}
\label{l:M=W/W_Q}
Let $G=\mathbf{G}(k)$ be the $k$-points of a simple algebraic group defined 
over a local field $k$ and $B=G/P$ be as in Lemma~\ref{L:semi-simple}.
Then in the sense of the Galois correspondence between $\mathcal{SG}(W_{G,B})$
and $\mathcal{Q}(B)$ the bijection between the closed objects 
is the correspondence between the special (parabolic) subgroups $W_{Q}<W$, 
and the closed factors $G/Q$, where $Q$ is a parabolic subgroup:
$T(W_{Q})=G/Q$ and $S(G/Q)=W_{Q}$.
\end{lem}
\begin{proof}
By the fact that closed factors are $G$-factors of $G/P$ we get that they
are all of the form $G/Q$ where $Q$ is a standard parabolic subgroup.
It is then a standard exercise in Lie theory to prove the correspondence
$G/Q\leftrightarrow W_Q$.
\end{proof}
The classical Weyl group $W=N_{\mathbf{G}}(\mathbf{S})/Z_{\mathbf{G}}(\mathbf{S})$
of a semi-simple group over $k$ acts on its $k$-relative root system,
and can be viewed as a finite Coxeter group acting on a real 
vector space $V$ of dimension $\text{rk}_{k}(\mathbf{G})$.
The group is simple iff its Weyl group acts irreducibly on $V$
(\cite[IV.14]{Borel:LinAlgGrps:91}).
\begin{lem}\label{l:coxeter-faithful}
Let $W$ be a finite irreducible Coxeter group,
and let $W'$ be a proper special (parabolic) subgroup.
Then $W$ acts faithfully on $W/W'$.
\end{lem}

\begin{proof}
Let $W\to\GL(V)$ be the Cartan representation, i.e the representation that
realizes $W$ as a reflection group ($V$ is a real vector space).
Denote by $V^+$ the Weyl chamber in $V$.
The representation $W\to\GL(V)$ is faithful and irreducible.
Indeed, if $V=V_1\oplus V_2$ is a $W$-invariant decomposition, each
reflection in $W$ must be trivial on one of the factors, and a reflection on the other.
The reflections which are trivial on $V_1$ commutes with the reflections 
which are trivial on $V_2$, then
by the irreducibility of $W$ as a Coxeter group,
we may assume that one of the set is empty.
That is, all the reflections act trivially on, say, $V_1$.
This means that $V_1$ is in the kernels of all the reflection, which
is an absurd.
To see that the representation is faithful it is enough to recall that the 
points in $V^+$ have trivial stabilizers.
Recall that $W'$ is the pointwise stabilizer of a face
of $V^+$.
$W'$ is a proper subgroup, hence this face is not $\{0\}$.
Denote the kernel of the action of $W$ on $W/W'$ by $K$.
$V^K \supset V^{W'} \neq \{0\}$.
As $K$ is normal in $W$, $V^K$ is $W$-invariant,
hence, by irreducibility, it is $V$.
We get that $K$ is the representation's kernel, which is trivial.
\end{proof}

\medskip

It is probably worth illustrating the concepts by a concrete example.
\begin{exmp}[$G=\SL_{3}$]
\label{EX:S3}
Let $k$ be a local field, $G=\SL_{3}(k)$ and $B=G/P$ where $P$ is the subgroup
of upper triangular matrices, i.e. the stabilizer of the line $\Span(e_{1})$ and the plane
$\Span(e_{1},e_{2})$ in the linear $G$-action on $V=k^{3}$. 
Then $B$ is the space of all pairs $(\ell,\Pi)$, a.k.a. \emph{flags},
of nested proper vector subspaces $\ell\subsetneq \Pi$ of $V$, i.e., configurations 
$\ell\subset \Pi$ of a line $\ell$ and a containing plane $\Pi$ in $V$.
Denoting by $A$ the diagonal subgroup of $G$, the $G$-space $G/A$ can be identified
with the space of all triples $(\ell_{1},\ell_{2},\ell_{3})$ of non-coplanar lines in $V$. 
Then the map $G/A\to B\times B$ given by 
\[
	(\ell_{1},\ell_{2},\ell_{3})\mapsto 
	((\ell_{1},\Pi_{1}=\ell_{1}\oplus\ell_{3}),(\ell_{2},\Pi_{2}=\ell_{2}\oplus\ell_{3}))
\]
is a $G$-equivariant Borel bijection onto a subset of full measure in $B\times B$
(cf. \ref{L:BxB=G/Z}).
Note that the symmetric group $S_{3}$, acting by permutations on the space of triples 
of such lines, is precisely the Weyl group in our sense (cf. \ref{P:W4ss}).
Let $\Phi:B\to C$ be a measurable quotient (such as a $\Gamma$-equivariant boundary map
$B\to \Prob(\Circle)$ arising from an action $\Gamma\to\Homeo_{+}(\Circle)$). 
Then $W'=S(\Phi)<W=S_{3}$ consists of those permutations $\sigma\in S_{3}$ for which
\[
	\Phi(\ell_{1},\ell_{1}\oplus\ell_{2})=\Phi(\ell_{\sigma(1)},\ell_{\sigma(1)}\oplus\ell_{\sigma(2)}).
\]
for a.e. $(\ell_{1},\ell_{2},\ell_{3})$.
For a fixed $W'<S_{3}$ there is a maximal quotient $\ol\Phi:B\to \ol{C}$ for which
the above holds, and therefore $\Phi$ factors through $\ol{\Phi}$.
These quotients are automatically $G$-equivariant\footnote{If $G$ has higher rank
and $\Phi$ is $\Gamma$-equivariant, then by Margulis Factor Theorem
$\Phi$ is $G$-equivariant, and so $C=G/Q$. However, we do not need to use 
this deep result, apealing instead to the Galois correspondence, which shows that
$\Phi$ always factors through $\ol{\Phi}:B\to\ol{C}$ of the form $\ol{C}=G/Q$ with 
$W_{Q}=\ol{W'}\supset W'$.} 
and the two non-trivial ones are the projections:
$\Phi_{1}(\ell,\Pi)=\ell$ and $\Phi_{2}(\ell,\Pi)=\Pi$.
For example if $S(\Phi)$ fixes $\ell_{1}$ interchanging $\ell_{2}$ with $\ell_{3}$, 
then $\Phi(\ell,\Pi)$ depends only on $\ell$
and $\ol{\Phi}=\Phi_{1}$; if $S(\Phi)$ interchanges $\ell_{1}$ with $\ell_{2}$ leaving $\ell_{3}$ fixed,
then $\Phi(\ell,\Pi)$ depends only on $\Pi$ and $\ol{\Phi}=\Phi_{2}$. In all other cases 
with a non-trivial $S(\Phi)$ the closure $\ol{\Phi}$ is the trivial map of $B$ to a point, 
which means that the original $\Phi$ is a.e. a constant map.
\end{exmp}

\medskip

In the context of $\tilde{A}_{2}$-groups we announce the following analogue 
of Proposition~\ref{P:key4simple}:
\begin{prop}
\label{P:key4A2}
Let $\Gamma$ be an $\tilde{A}_{2}$ group, $(B,\nu)$ a $\Gamma$-boundary 
as in \ref{P:A2bnd}, and $W<W_{\Gamma,B}$ a group isomorphic to $S_{3}$ 
as in  \ref{P:W4A2}.
Then for any quotient map $\Phi:B\to C$ the $W$-action on $W/S(\Phi)$ is faithful,
unless $\Phi$ is the constant map to the singleton.
\end{prop}
%

\section{Circle peculiarities} \label{S:Circle}
In this section we collect some special properties of the circle.
We think of the circle as a compact 1-manifold, or as $\RR/\ZZ$,
but it can be thought of as a completion of \emph{the unique} 
countable set with a dense cyclic order.
This point of view will come handy in Section~\ref{S:Appendix}. 

\subsection{Dynamics on the circle} \label{circle-dynamics}

Let $\ro:\Gamma\to \Homeo_{+}(\Circle)$ be a homomorphism, i.e. an action
of some group $\Gamma$ on the circle $\Circle$ by orientation preserving
homeomorphisms. We list some classical facts that will be used below.

\medskip

\begin{prop}[Minimal Sets]
\label{P:minimal-sets}
Suppose $\Gamma$ has only infinite orbits.
Then either 
\begin{enumerate}
\item		
	$\Gamma$ acts minimally on the whole circle $\Circle$, or 
\item
	$\Gamma$ has a unique minimal set $K\subset\Circle$,
	in which case $K$ is homeomorphic to a Cantor set. 
\end{enumerate}
\end{prop}

\begin{lem}[Semi-conjugacy]
\label{L:semi-conj}
Let $K\subset\Circle$ be a closed $\Gamma$-invariant set. 
Then there exists a degree one continuous monotonic map $f:\Circle\to\Circle$
and $\ro^{\prime}:\Gamma\to \Homeo_{+}(\Circle)$ so that 
\[
	f(K)=\Circle\qquad\text{and}\qquad
	\ro^{\prime}(\gamma)\circ f=f\circ \ro(\gamma),\qquad (\gamma\in\Gamma).
\]
If $\Gamma$ has no fixed points in $K$ then $\Ker(\ro)=\Ker(\ro^{\prime})$.
\end{lem}
%


\begin{lem}[Invariant measure]
\label{L:inv-meas}
If a subgroup $\Gamma<\Homeo_{+}(\Circle)$ leaves invariant a probability measure
$\mu\in\Prob(\Circle)$, then either $\Gamma$ has a finite orbit on $\Circle$,
or $\Gamma$ is is semi-conjugate to a dense subgroup of the rotation group $\SO_{2}(\RR)$.
In the latter case the semi-conjugacy can be replaced by a conjugacy if
$\mu$ has full support.
\end{lem}
\begin{proof}
Let $\mu$ has atoms, then its purely atomic part $\mu_{a}$ is $\Gamma$-invariant.
If it is non-empty, then the set of $\mu$-atoms of maximal weight is a finite non-empty
$\Gamma$-invariant set. Hence $\Gamma$ has at least some finite orbits.

Assuming $\Gamma$ has only infinite orbits, $\mu$ has no atoms.
Viewing $\Circle$ as $\RR/\ZZ$ define the function $f:[0,1]\to[0,1]$ by $f(t)=\mu([0,t])$.
It is a continuos, monotonically increasing function with $f(0)=0$, $f(1)=1$.
So it can be viewed as a continuous degree one map $f:\Circle\to\Circle$. 
The assumption that $\Gamma$ preserves $\mu$ yields that $f$ is a semi-conjugacy
of $\Gamma$ into the group of rotations. 
If $\supp(\mu)=\Circle$ then $f$ is a homeomoprhism,
and hence $\Gamma$ is \emph{conjugate} in 
$\Homeo_{+}(\Circle)$ into $\SO_{2}(\RR)$.

\end{proof}

\medskip

The above stated facts are classical and were probably known to Poincar\'e.
The following useful result is more recent, it is based on an argument of 
Margulis \cite{Margulis:Tits:00} 
with an addition due to Ghys \cite[see pp.361--362]{Ghys:cicrcle:01}. 
For its formulation we need a term to describe a certain type of dynamics.
There are several commonly used properties of group actions on general (compact)
spaces, which amount to the same property when applied to the circle.
We have picked the term \emph{minimal and strong proximal} action:
\begin{defn}
\label{D:str-prox}
We shall say that an action  $\ro:\Gamma\to\Homeo_+(\Circle)$ of some group 
on the circle is \emph{minimal and strongly proximal} if given any proper closed 
subset  $K\subsetneq \Circle$ and any non-empty open set 
$U\subset \Circle$ there exists $\gamma\in\Gamma$
with $\ro(\gamma)K\subset U$.
(Of course, we may take $K$ and $U$ to be just proper closed and open arcs on the circle).
\end{defn}
\begin{thm}[Margulis, Ghys]
\label{T:Margulis-Ghys}
Let $\ro:\Gamma\to\Homeo_+(\Circle)$ be a minimal action of some group $\Gamma$.
Then either 
\begin{enumerate}
\item
	$\ro(\Gamma)$ is equicontinuous, in which case $\ro(\Gamma)$ can be conjugated into
	the group of rotations, or
\item
	the centralizer $C$ of $\ro(\Gamma)$ in $\Homeo_+(\Circle)$ is a finite cyclic group,
	and the quotient action
	\[
		\ro^\prime:\Gamma\overto{\ro}\Homeo_+(\Circle)
		\overto{}\Homeo_+(\Circle/C)
	\]
	is minimal and strongly proximal.
\end{enumerate}
\end{thm}
\begin{cor}
\label{C:no-inv-meas}
Any non-elementary continuous action of an arbitrary group on the circle  
is semi-conjugate to an action which is minimal and strongly proximal.
\end{cor}
\begin{proof}
Let $\ro:\Gamma\to\Homeo_{+}(\Circle)$ be an action without invariant probability measures.
In particular it has no finite orbits. This action is then semi-conjugate to a minimal one
using a degree one monotonic map $f:\Circle\to\Circle$. 
The resulting action $\ro^{\prime}$ cannot be conjugate into rotations because 
rotations leave the Lebesgue measure invariant, and its lift (well defined)
to the minimal set for $\ro$ would be $\ro(\Gamma)$-invariant.
Passing to a finite-to-one quotient of $\ro^{\prime}$ as in (2) of Theorem~\ref{T:Margulis-Ghys} 
we obtains a minimal strongly proximal action.
\end{proof}
\begin{proof}[Sketch of the proof of Theorem~\ref{T:Margulis-Ghys}]
Say that an arc $I\subset\Circle$ is contractible if there exists a sequence $\gamma_n\in\Gamma$
so that the length of $\ro(\gamma_n)(I)$ goes to $0$. 

Alternative (1): there are no non-trivial contractible arcs.
It is not difficult to see that in this case $\ro(\Gamma)$ is equicontinuous.
Therefore the closure $Q=\overline{\ro(\Gamma)}$ in $\Homeo_{+}(\Circle)$ is a compact group
by Arzela-Ascoli.
The pushforward of its Haar measure is a $\ro(\Gamma)$-invariant probability measure. 
It has full support on the circle because the action is assumed to be minimal.
Lemma~\ref{L:inv-meas} completes the argument. 

Alternative (2): there exist non-trivial contractible arcs. 
Using minimality one shows that contractible intervals cover the whole circle.
Denoting by $\pi:\RR\to\RR/\ZZ=\Circle$ the projection, define a function 
$\tilde\theta:\RR\to\RR$ by 
\[
	\tilde\theta(x)=\sup\{ y \mid \pi([x,y]) \text{ is contractible} \}
\]
Note that $\tilde{\theta}$ is well defined, satisfies $x<\tilde{\theta}(x)\le x+1$ and is non-decreasing.
The case of $\ro(\Gamma)$ being strongly proximal corresponds precisely to $\tilde\theta(x)=x+1$.
Using minimality it is shown that $\tilde\theta$ has no jumps, and no flat intervals. 
Clearly $\tilde\theta$ commutes with $\ZZ$, and so defines $\theta\in\Homeo_{+}(\Circle)$;
the latter commutes with $\ro(\Gamma)$. 
Let $C$ denote the cyclic group $C$ generated by $\theta$.
As $\theta$  cannot have irrational rotation number
(for then $\Gamma$ would be equicontinuous), $C$ is finite.
The quotient by $C$ defines a minimal and strongly proximal action.
\end{proof}

\subsection{Boundary maps for Circle actions}

We start by a small digression describing a natural embedding
of the space $\Meas_{c}(\Circle)$ of all finite, signed, \emph{continuous} measures on the circle 
into $C(\Sphere)$ -- the space of continuous functions on the sphere.

We first observe that the space of all closed proper 
arcs (arcs with two different end points)
on the circle $\Circle$,
endowed with the Hausdorff metric is
homeomorphic to 
$\Circle\times\Circle\setminus\Delta(\Circle)$,
which is a cylinder $\Circle\times (0,1)$.
The action of $\Homeo_+(\Circle)$ on the circle extends to a continuous
action on the cylinder.
Any non-atomic measure $\mu$ on $\Circle$, defines a function 
on the cylinder, by an evaluation of the $\mu$-measure of arcs.
This function, denoted $f_{\mu}$, is continuous, vanishes at one side $\Circle\times(0,\epsilon)$,
and tends to $\mu(\Circle)$ on the other side $\Circle\times (1-\epsilon,1)$,
because $\mu$ is assumed to have no atoms. 
Hence, $f_{\mu}$ is defined on the two-point compactification of the cylinder,
namely on $\Sphere$.
This gives a $\Homeo_+(\Circle)$ equivariant embedding 
\[
	\Meas(\Circle)\overto{i} C(\Sphere)
\]
which is a positive linear operator of norm $1$ with respect to the variation 
norm on $\Meas(\Circle)$
and $\max$-norm on $C(\Sphere)$. 
This immediately gives:
\begin{lem}
\label{L:dist}
The pullback of the $\max$-norm on $C(\Sphere)$ 
\[
	d(\mu_1,\mu_2)=\max_{p\in\Sphere} |f_{\mu_1}(p)-f_{\mu_2}(p)|
		=\max_{I\subset \Circle}|\mu_1(I)-\mu_2(I)|
\]
is a $\Homeo_{+}(\Circle)$ invariant metric on the space $\Meas_c(\Circle)$ 
of all continuous measures on the circle. 
\end{lem}

The following fact makes circle actions easy to analyze using boundary maps:
\begin{thm}[Reduction to atomic measures]  
\label{T:reduction2atoms} 
Let $\ro:\Gamma\to\Homeo_+(\Circle)$ be a minimal and strongly proximal action
of some group $\Gamma$ on the circle, 
and let $(B,\nu)$ be some doubly ergodic measurable $\Gamma$-space.
Assume that there exists a measurable $\Gamma$-equivariant map 
$\Phi:B\to\Prob(\Circle)$. 
Then
\begin{enumerate}
\item
	Such $\Phi$ is unique (as a map defined up to $\nu$-null sets).
\item
	There exists $k\in\NN$ such that for $\nu$-a.e. $x\in B$ the measure $\Phi(x)$ is atomic, 
	supported and equidistributed on a $k$-point set $A_{x}\subset \Circle$.
\item
	For a.e. $x,y\in B\times B$ the $k$-point sets $A_x$ and $A_y$ 
	are "unlinked", i.e., belong to complimentary arcs of the circle.
\item
	In fact, the value of $k$ in (2) is either $k=1$ or $k=2$.
\end{enumerate}
\end{thm}
\begin{rem}
We do not know any example where $k=2$. 
In fact, it can be shown that $k=1$ in many situations: 
it was proved for the case of $(B,\nu)$ being the Poisson boundary of a generating random walk on 
$\Gamma$ by Antonov \cite{Antonov} (see \cite[Proposition 5.7]{DeroinKleptsynNavas:05}), 
and independently by Furstenberg (unpublished); 
it can be shown under a more general assumption that $(B,\nu)$ is an SAT space, 
or if there is an (infinite) invariant measure in the class of $\nu\times\nu$ on $B\times B$.
If indeed $k=1$, then the statement of the Theorem is simply: \emph{there exists a measurable 
$\Gamma$-equivariant map $\phi:B\to\Circle$, and $\Phi(x)=\delta_{\phi(x)}$ 
is the unique $\Gamma$-equivariant map $B\to\Prob(\Circle)$}. 
Such a statement could be convenient, but for the purposes of this paper we shall be satisfied with properties (1)--(3) as stated.  
\end{rem}

For an integer $k$ let $\Circle_{k}$ denote the collection of all $k$-point subsets of the circle. 
We start by some general observations about the behavior of measurable $\Gamma$-equivariant maps
$B\to\Circle_{k}$.
\begin{lem}
\label{L:squeez}
Let $A:B\to\Circle_{k}$, $x\mapsto A_{x}$ be a measurable $\Gamma$-equivariant map.
Then for any non-empty open set $U\subset\Circle$ the set of $x\in B$ with
$A_{x}\subset U$ has positive $\nu$-measure.
\end{lem}
\begin{proof}
For any fixed point $p\in\Circle$ we have $\nu\{x\in B\mid p\in A_{x}\}<1$, because otherwise 
$p$ would be fixed by a finite index subgroup of $\ro(\Gamma)$. 
Then for a sufficiently large proper closed arc $K\subset \Circle\setminus\{p\}$
the set $E=\{x\in B\mid A_{x}\subset K\}$ has $\nu(E)>0$.
Due to  strong proximality and minimality of the action there exists $\gamma\in\Gamma$ with 
$\ro(\gamma)K\subset U$. Then $F=\gamma E$ has $\nu(F)>0$ and for $x\in F$
\[
	A_{x}\subset \gamma K\subset U.
\]
\end{proof}
\begin{lem}
\label{L:uniq-Ax}
Let $A:B\to \Circle_{k}$ and $A^{\prime}:B\to \Circle_{k^{\prime}}$ be measurable 
$\Gamma$-equivariant maps. 
Then $k=k^{\prime}$ and $A_{x}=A^{\prime}_{x}$ for $\nu$-a.e. $x\in B$.
\end{lem}
\begin{proof}
We need some notations.  
Given a pair $(S,S^{\prime})\in\Circle_{k}\times\Circle_{k^{\prime}}$
mark each of the points $p\in S\cup S^{\prime}$ by one of the symbols $\{+,-,*\}$ 
according to whether $p$ lies in
\[
	S\setminus S^{\prime}, \qquad S^{\prime}\setminus S,\qquad\text{or}\qquad
	S\cap S^{\prime}.
\]
The set $S\cup S^{\prime}$ divides the circle into $\le k+k^{\prime}$
open arcs $(p,q)$ with $p,q\in S\cup S^{\prime}$ (the interval is named according to the fixed 
orientation of the circle).
There are $\{+,-,*\}^{2}$ possible markings of the endpoints of such arcs. 
We shall divide the arcs into \emph{positive}, \emph{negative}
and neutral, by grouping types $(-,+)$, $(-,*)$, $(*,+)$ into the positive class,
and $(+,-)$, $(*,-)$, $(+,*)$ into the negative class, with the rest being neutral.
The point of these definitions is that for any configuration of $S,S^{\prime}$ with $S\neq S^{\prime}$
there will be at least one positive arc and at least one negative one.
The markings of the points and the labeling of the arcs are clearly
$\Homeo_{+}(\Circle)$-invariant.

Now consider the measurable map 
\[
	(x,y)\in B\times B\ \mapsto\ (A_{x}, A_{y}, A^{\prime}_{x})\in 
	\Circle_{k}\times \Circle_{k}\times \Circle_{k^{\prime}}.
\]
The pair $(A_{x},A^{\prime}_{x})$ divides the circle into positive, 
negative and neutral arcs as above. 
Let $P_{x}$ denote the union of the positive arcs and $N_{x}$
the union of the negative ones. These are measurable $\Gamma$-equivarinat
families of open sets, which are non-empty iff $A_{x}\neq A^{\prime}_{x}$.

Define $f:B\times B\to \{-1,0,1\}$ by letting $f(x,y)=1$ if $A_{y}\subset P_{x}$,
$f(x,y)=-1$ if $A_{y}\subset N_{x}$; and $f(x,y)=0$ in all other cases.
This function is clearly measurable and $\Gamma$-invariant,
hence $\nu\times\nu$-a.e. constant. 
We claim that this is possible only if $A_{x}=A^{\prime}_{x}$ a.e.
in which case $P_{x}=N_{x}=\emptyset$ and $f=0$.
Otherwise for a positive measure set $B_{0}$ of $x\in B$ the sets 
$P_{x}$ and $N_{x}$ are non-empty and by Lemma~\ref{L:squeez} the sets
\[
	E_{x}=\{y\in B \mid A_{y}\subset P_{x}\},\qquad F_{x}=\{y\in B \mid A_{y}\subset N_{x}\}
\]
have $\nu(E_{x})>0$ and $\nu(F_{x})>0$. By Fubini, this means that 
$f$ takes both $+1$ and $-1$ values on sets of positive $\nu$ measure, 
which is a contradiction.
\end{proof}

\begin{proof}[Proof of Theorem~\ref{T:reduction2atoms}]
Let $\Phi:B\to\Prob(\Circle)$ be a given measurable equivariant map.
We recall the elegant argument of Ghys showing that $\Phi$ is necessarily atomic.

The decomposition $\mu=\mu_{a}+\mu_{c}$ of measures
$\mu\in\Prob(\Circle)$ into the purely atomic and the continuous
(= the complement of the atomic) parts, is clearly $\Homeo_{+}(\Circle)$-invariant.
So their corresponding weights are $\nu$-a.e. constant, and the claim
is that $\Phi(x)$ has no continuous part.

Had it been different, the normalized continuous part $\Phi_c(x)$ of $\Phi(x)$
would give a $\Gamma$-equivariant measurable map $B\to\Prob_c(\Circle)$.
Let $d$ be the $\Homeo(\Circle)$-invariant metric as in Lemma~\ref{L:dist}.
Then $d(\Phi_c(x),\Phi_c(y))$ is a measurable $\Gamma$-invariant function on $B\times B$
and hence is $\nu\times\nu$-a.e. a constant $d_0$.
This constant has to be $d_0=0$. 
Indeed otherwise we could cover 
$\Prob_c(\Circle)$, imbedded in $C(\Sphere)$ via $\mu\mapsto f_\mu$, by countably 
many sets $\{U_{n}\}$ of diameter $<d_{0}$ and define
\[
	E_{n}=\{x\in B \mid f_{\Phi(x)}\in U_{n}\}.
\]
Then up to $\nu$-null set $B$ is covered by $\bigcup_{1}^{\infty} E_n$,
while $\nu(E_n)=0$ because for $\nu$-a.e. $x\in E_{n}$ for $\nu$-a.e. $y\not\in E_{n}$.
The contradiction yields $d_0=0$. Thus $\Phi(x)$ is a fixed measure $\mu_0\in\Prob_c(\Circle)$.
This continuous measure is therefore $\ro(\Gamma)$-invariant.
This argument, originally due to Ghys \cite[207-210]{Ghys:lattices:99}, applies to
arbitrary actions. 
However our assumption of strong proximality rules out existence of a 
$\ro(\Gamma)$-invariant  measure. 
Thus $\Phi(x)$ is a purely atomic measure for $\nu$-a.e. $x\in B$.

For any $\epsilon>0$ the set $A_{x,\epsilon}$ of all the atoms of the atomic probability
measure $\Phi(x)$ with weight $\ge \epsilon$ is finite (bounded by $\epsilon^{-1}$).
We also have 
\[
	A_{\gamma.x,\epsilon}=\ro(\gamma)A_{x,\epsilon}.
\] 
The cardinality $|A_{x,\epsilon}|$ of these sets is a $\Gamma$-invariant function, and therefore
a.e. constant $k_{\epsilon}$. 
There is $\epsilon_{0}>0$ so that $k_{\epsilon}\ge1$ for $\epsilon\in(0,\epsilon_{0}]$.
Lemma~\ref{L:uniq-Ax} shows that the set $A_{x,\epsilon}$ ($0<\epsilon\le\epsilon_{0}$)
depends only on $x$ and can be denoted $A_{x}$. 
It follows that $\epsilon_{0}=k^{-1}$ where $k=|A_{x}|$ and
\[
	\Phi(x)=\frac{1}{k} \sum_{a\in A_{x}}\delta_{a}.
\]
This proves assertions (1) and (2).

To show (3) consider the set $E$ of all pairs $(x,y)\in B\times B$ for which
the sets $A_{x}$ and $A_{y}$ can be separated on the circle. 
For $\nu$-a.e. $x$ the fiber $E_{x}=\{y\in B \mid (x,y)\in E \}$ has positive measure,
using Lemma~\ref{L:squeez}. Thus $\nu\times\nu (E)>0$. 
But $E$ is $\Gamma$-invariant, and therefore is conull.  

(4). Note that for $k\ge 3$ any family $\sA$ of unlinked $k$-point sets on the circle
is at most countable. To see this, define  \emph{separation} of a $k$-point set $A\subset\Circle$ to be
\[
	{\rm sep}(A) =\min\{ d(a,a^{\prime}) \mid a\neq a^{\prime}\in A\}.
\]
It suffices to show that for any $\epsilon>0$ the collection $\sA_{\epsilon}=\setdef{A\in\sA}{{\rm sep}(A)\ge\epsilon}$ is finite. This is indeed so, in fact the cardinality of $\sA_{\epsilon}$ is bounded by $(k-2)/\epsilon$, because the union of all sets in $\sA_{\epsilon}$ divides the circle into arcs, 
with each set $A\in\sA_{\epsilon}$ 
contributing at least $(k-2)$ arcs of length $\ge \epsilon$.
\end{proof}
%

\subsection{The orbit structure of  $(\Circle)^{n}$ and $(\Circle_k)^m$}


Fix a natural number $n$.
The space $(\Circle)^{n}$
is acted uppon by a diagonal action of $\Homeo_{+}(\Circle)$,
and by coordinate permutation action of the symmetric group $S_n$.
These actions commute.

Let $d:(\Circle)^{n} \to \NN$ be the map
$$(x_1,\ldots,x_n)\mapsto |\{x_1,\ldots,x_n\}| $$
This map is $\Homeo_{+}(\Circle) \times S_n$ invariant.
It is not hard to see that the $\Homeo_{+}(\Circle)$ orbit of $x \in (\Circle)^{n}$ is a manifold of dimension
$d(x)$, and that
the closure of an $\Homeo_{+}(\Circle)$ orbit consists of orbits of lower dimensions.
In particular there are $n!$ open orbits, and a unique closed orbit --
the diagonal.

Assume from now on that $n=km$.
Let $L=S_k\times S_k\times \cdots\times S_k < S_{km}$.
Denote 
\[ X=\{x\in (\Circle)^{n}~|~\Stab_{S_n}(x)\cap L=\{1\}\} \]
We identify $X/L$ with $(\Circle_k)^m$.
The map $d$ decends to a map, denoted $d$ as well, $d:X\to\NN$.
We readily get the following lemma.

\begin{lem} \label{l:tame}
The $\Homeo_{+}(\Circle)$ orbit of every point in $(\Circle_k)^{m}$ is a manifold.
The closure of each orbit consists of orbits of lower dimensions.
In particular the quotient topology on $(\Circle_k)^{m}/\Homeo_{+}(\Circle)$ is $T_0$,
and the $\Homeo_{+}(\Circle)$ action on $(\Circle_k)^{m}$ is tame.
\end{lem}

\begin{lem} \label{l:aut(orbit)}
Let $\mathcal{O}$ be an orbit of $H=\Homeo_{+}(\Circle)$  in
$(\Circle_{k})^{m}$.
Then the group $\Aut_H(\mathcal{O})$ of permutations of $\mathcal{O}$
commuting with $H$, is a finite cyclic group.
\end{lem}
\begin{proof}
Fix a point 
$\vec{A}=(A_{1},\dots,A_{m})\in \mathcal{O}$, $A_{i}\in \Circle_{k}$ and
let $A=A_{1}\cup\dots\cup A_{m}$.
Denote $H_{\vec{A}}=\Stab_H(\vec{A})$, $F_A=\Fix_H(A)$, and $H_A=\Stab_H(A)$.
Obviously $ F_A < H_{\vec{A}} < H_A $,
and $F_A$ is normal in $H_A$, with a finite cyclic quotient.
Observe that every finite set of the circle which is stabilized by $F_A$ is already in $A$.
It follows that for every intermidate subgroup $F_A<S<H_A$,
$N_H(S)$ is contained in $H_A$ as well.
Furthemore, $N_H(S)/S$ is a finite cyclic group, as a subquotient of $H_A/F_A$.
In particular, $N_H(H_{\vec{A}})/H_{\vec{A}}$ is finite and cyclic.
(In fact, $F_A$ can be seen to be characteristic in $H_A$ -- 
it is the kernel of its profinite "completion" -- but we will not use this fact here).

By identifying $\mathcal{O}$ with $H/H_{\vec{A}}$, 
and using Lemma~\ref{L:Aut=N/H},
we have $\Aut_H(\mathcal{O}) \cong N_H(H_{\vec{A}})/H_{\vec{A}}$, 
thus the proof is complete.
\end{proof}

\subsection{Weyl group and the circle}
\label{SS:cyclicity}
Recall the definitions and the construction of an abstract Weyl group in \ref{SS:Weyl}. 
\begin{thm} 
\label{T:cyclic}
Let $G$ be a general \lcsc group, $\Gamma<G$ a lattice with 
a minimal and strongly proximal action $\ro:\Gamma\to\Homeo_{+}(\Circle)$ on a circle.
Let $(B,\nu)$ be a $G$-boundary and $\Phi:B\to\Prob(\Circle)$
a $\Gamma$-equivariant map.
Let $W$ be a finite subgroup of $W_{G,B}$, containing $w_0$.
Then $W$ acts on $W/S(\Phi)$ through a non-trivial finite cyclic quotient group.
\end{thm}
\begin{rem}
The assumption that $\ro(\Gamma)$ is minimal and strongly proximal
is superfluous. But we shall use the result only in the simpler context 
of a minimal strongly proximal case.
\end{rem}
The theorem is an easy consequence of the constructions in \ref{S:Galois}
combined with Lemmas~\ref{l:tame} and ~\ref{l:aut(orbit)}.

\begin{proof}
From Theorem~\ref{T:reduction2atoms} we know
that $\alpha$ takes values in atomic measures and may be viewed as a 
measurable $\Gamma$-equivariant map
\[
	\Phi:B\to\Circle_k.
\] 
Consider the corresponding $\Gamma\times W$-equivariant map
\[
	\vec{\Psi}: B\times B \overto{} \prod_{W/S(\Phi)} \Circle_{k} 
	= (\Circle_{k})^{m}
\]
given by Construction \ref{cons:vec-map}, where we denote $|W|=m$.
By Proposition~\ref{prop:boundary-lattice}, $\Gamma$ acts ergodically on $B\times B$,
hence the pushforward measure class $\vec{\Psi}_{*}([\nu\times\nu])$ is an ergodic 
measure class on $(\Circle_{k})^{m}$,
with respect to $\ro(\Gamma)<H$.
As the action of $H$ on $(\Circle_{k})^{m}$ is tame (Lemma~\ref{l:tame}),
we get that $\vec{\Psi}_{*}([\nu\times\nu])$ is supported on a single $H$-orbit, 
say $\mathcal{O}$.
Since $W$-action on $B\times B$ preserves the class $[\nu\times\nu]$,
the $W$-action on $\prod_{W/S(\Phi)} \Circle_{k}$ preserves
the measure class $\vec{\Psi}_{*}([\nu\times\nu])$.
It follows that $\mathcal{O}$ is $W$-stable.
By Lemma~\ref{l:aut(orbit)} the image of $W\to \Aut_{H}(\mathcal{O})$
is a finite cyclic group.
The $W$-action on $\mathcal{O}$ is given by permutation of coordinates,
and is obtained from its action on $W/S(\Phi)$.
We are just left to remark that since $\ro(\Gamma)$ cannot have 
an invariant measure in $\Prob(\Circle)$,
$\Phi$ is a non-trivial factor map,
hence by Lemma~\ref{l:long-element},
$w_0 \notin S(\Phi)$, so the $W$ action on $W/S(\Phi)$
is not trivial.
\end{proof}

\medskip

\subsection{Extending actions from a dense subgroup}
\begin{lem}[Extending homomorphisms]
\label{L:extending}
Let $\Lambda$ be a dense subgroup in some \lcsc group $G$,
and $\ro:\Lambda\to\Homeo_{+}(\Circle)$ be a minimal and strongly proximal 
action (see \ref{D:str-prox}).
Suppose that there exists a measurable $G$-space $(B,\nu)$
which admits a measurable $\Lambda$-equivariant map
$\Phi:B\to\Prob(\Circle)$. 

Then $\ro$ extends (uniquely) to a continuous homomorphism 
$\bar{\ro}:G\overto{}\Homeo_{+}(\Circle)$ with a \lcsc image
$\bar{\ro}(G)<\Homeo_+(\Circle)$. 
\end{lem}

\begin{proof}
It suffices to prove that $\ro$ is continuous with respect to the topology
on $\Lambda$ induced from $G$, that is to prove that for any sequence 
$\lambda_{n}\in\Lambda$ with $\lambda_{n}\to e$ in $G$, we have
$\ro(\lambda_n)\to \Id$ in $\Homeo_+(\Circle)$.
Indeed this would imply that for any $g\in G$ and any sequence $\lambda_n\to g$
the corresponding sequence $\{\ro(\lambda_n)\}$ is Cauchy, and therefore
uniquely defines an element $\bar{\ro}(g)$ in the Polish group $\Homeo_+(\Circle)$.
This also shows that $g\mapsto \bar{\ro}(g)$ is well defined, and is a continuous group
homomorphism. 
The image of a \lcsc group under a continuous homomorphism is easily seen to be 
\lcsc as well.

So fix some sequence $\lambda_{n}\in\Lambda$ with $\lambda_{n}\to e$ in $G$.
Upon passing to a subsequence we may assume that $\ro(\lambda_{n})$ converges
point-wise to a monotonic, but not necessarily continuous, map $f:\Circle\to\Circle$
(choose a countable dense set of points $\{p_{i}\}$, 
and applying the diagonal argument to $\ro(\lambda_{n}).p_{i}$ 
to extract subsequence converging for each $i$, 
let $f(p_{i})=\lim\ro(\lambda_{n}).p_{i}$, and extend by monotonicity).
If $f(p)=p$ for all $p$ then $\ro(\lambda_{n})$ converge (uniformly) to the identity.

Otherwise the set $\{p\in\Circle \mid f(p)\neq p\}$ is uncountable and contains 
a point of continuity for $f$, and therefore we can find non-trivial disjoint open 
arcs $I,J\subset\Circle$ with $f(I)\subset J$. Let 
\[
	E=\{x \in B \mid \supp(\Phi(x))\subset I\},\qquad
	F=\{x \in B \mid \supp(\Phi(x))\subset J\}.
\]
These are disjoint sets of positive $\nu$-measure (Lemma~\ref{L:squeez}).
For $\nu$-a.e. $x\in E$ we have $\lambda_{n}.x\in F$ starting from some $n=n(x)$.
At the same time
\[
	\nu(\lambda_{n}E\triangle E)\to 0\qquad\text{as}\qquad n\to\infty
\]
because $\lambda_{n}\to e$ in $G$. 
This holds for the same reasons as (and easily follows from) the continuity
of the quasi-regular $G$-representation on, say $L^{2}(B,\nu)$.
The contradiction completes the argument.
\end{proof}

\subsection{Locally compact subgroups acting on the circle}
$\Homeo_{+}(\Circle)$ is a Polish group, that is a complete separable topological group.
Understanding discrete subgroups of $\Homeo_{+}(\Circle)$ is a major open problem.
However non-discrete \lcsc subgroups of $\Homeo_{+}(\Circle)$ are well understood 
using Montgomery-Zippin theorem. Any non-discrete \lcsc subgroup of 
$\Homeo_{+}(\Circle)$ which acts minimally is conjugate to either $\SO(2)$ -- the group
of rotations, or to $\PSL_{2}(\RR)_{k}$ -- a $\ZZ/k\ZZ$ central extension of $\PSL_{2}(\RR)$
corresponding to the lifting of the standard $\PSL_{2}(\RR)$ action to the $k$-fold covering 
of the circle.
For minimal and strongly proximal actions we therefore get the following Theorem
(see \cite[pp. 345, 348]{Ghys:cicrcle:01} or \cite[pp. 51--53]{Furman:MM:01} for details).
\begin{thm}
\label{T:lcsubgroups}
Let $L<\Homeo_{+}(\Circle)$ be a non-discrete \lcsc subgroup which acts 
minimally and strongly proximally on the circle.
Then $L$ is conjugate to $\PSL_{2}(\RR)$ in $\Homeo_{+}(\Circle)$.
\end{thm}

\section{Higher rank Superrigity}
\label{S:mainproofs}
\subsection{Proof of Theorem \ref{T:product}}

\begin{proof}
Let $\Gamma$ be an irreducible lattice in the product 
$G=G_{1}\times \cdots\times G_{n}$ 
of $n\ge 2$  \lcsc  groups, and let  $\ro:\Gamma\to \Homeo_{+}(\Circle)$
be a homomorphism so that $\ro(\Gamma)$ is not elementary.
Thus, in particular, $\Gamma$ is not amenable, and so is $G$, and hence 
at least one of the its factors $G_{i}$.
After reordering we may assume that for some $1 \le m\le n$ the factors 
$G_1,\ldots,G_m$ are non-amenable, while $G_{m+1},\ldots,G_n$ are amenable. 
Denote $G'=G_1 \times \cdots \times G_m$ (it is possible that $m=n$ and $G'=G$).
For every $i=1,\ldots,m$ let $B_i$ be a $G_i$-boundary, and set 
$B=B_1\times\cdots\times B_m$.
By Proposition~\ref{P:products-of-boundaries} $B$ is a $G'$-boundary as 
well as a $G$-boundary.
Recall the notation introduced in observation~\ref{O:Weyl4products},
and let $W=\langle \tau_1,\ldots,\tau_m\rangle \simeq (\ZZ/2\ZZ)^m$ be the 
corresponding subgroup of $W_{G,B}$.

By Corollary~\ref{C:no-inv-meas}, $\ro:\Gamma\to\Homeo_{+}(\Circle)$ is semi-conjugate to 
a minimal and strongly proximal action $\ro_{1}:\Gamma\to \Homeo_{+}(\Circle)$. 
The $G$-boundary $(B,\nu)$ is also a $\Gamma$-boundary 
(see Proposition~\ref{prop:boundary-lattice}). 
Amenability of the $\Gamma$-action on $(B,\nu)$  
yields a $\Gamma$-equivariant map $\Phi:B\to \Prob(\Circle)$.
Theorem~\ref{T:cyclic} tells us that the action of $W$ on $W/S(\Phi)$ 
factors through a non-trivial finite \emph{cyclic} quotient group.
Lemma~\ref{L:WPhi4products} then implies that 
$\Phi:B\to  \Prob(\Circle)$ factors through some $B_{i}\to \Prob(\Circle)$.
By the irreducibility assumption on $\Gamma$, its projection to $G_i$ is dense.
Hence, we are in a position to apply Lemma~\ref{L:extending}.
It follows that $\ro_{1}$ can be extended to a continuous homomorphism 
$\bar{\ro}_{1}:G_i\overto{}\Homeo_{+}(\Circle)$.
An image of a \lcsc group in a Polish group under a continuous homomorphism is \lcsc
Since $L=\bar{\ro}_{1}(G_i)<\Homeo_+(\Circle)$ contains $\ro_{1}(\Gamma)$, it 
acts minimally and strongly proximally. 
The assumption that $G_{i}$ have no infinite discrete quotients, implies that
$L$ is conjugate to $\PSL_{2}(\RR)$ (Theorem~\ref{T:lcsubgroups}), 
and $\bar{\ro}_{1}$ is conjugate to an epimorphism 
$p:G_{i}\to\PSL_{2}(\RR)$. 
Finally, the orginal action $\ro:\Gamma\to\Homeo_{+}(\Circle)$ 
is semi-conjugate to the action 
\[
	\ro_{0}:\Gamma< G\overto{\pi_{i}} G_{i}\overto{p} \PSL_{2}(\RR)<\Homeo_{+}(\Circle)
\]
as claimed.
\end{proof}
\begin{rem}
\label{R:discreteimage}
Note that the assumption that $G_{i}$ have no non-trivial open normal subgroups,
is used only to identify $\bar\ro_{1}(G_{i})$ as $\PSL_{2}(\RR)$.
Without this assumption, we could still claim that a non-elementary 
$\ro:\Gamma\to\Homeo_{+}(\Circle)$ is semi-conjugate to $\ro_{1}$
which extends to $G$ and factors through some 
$G_{i}\overto{p}\Homeo_{+}(\Circle)$
with $p(G_{i})$ being either $\PSL_{2}(\RR)$ or \emph{discrete}.
\end{rem}

\medskip

\subsection{A simple proof of Theorems~\ref{T:higherrank} and  \ref{T:A2}}
\label{S:Ghys-A2}

%
\begin{proof}[Proof of theorem \ref{T:higherrank}]

Let $k$ be a local field, $\mathbf{G}$ be a $k$-connected simple algebraic 
group defined over $k$, ${\rm rk}_{k}(\mathbf{G})\ge 2$. 
Let $\Gamma$ be a lattice in the \lcsc group $G=\mathbf{G}(k)$, 
and $\ro:\Gamma\overto{}\Homeo_{+}(\Circle)$ a homomorphism.
If the action $\ro$ is elementary then by Lemma~\ref{L:inv-meas},
$\Gamma$ has a finite orbit on $\Circle$, as $\Gamma$ has finite Abelianization.
We will assume that $\ro$ is not elementary, and show that this leads to a contradiction.

By corollary~\ref{C:no-inv-meas}, $\ro$ is semi-conjugate to an action which
is minimal and strongly proximal, which we still denote by $\ro$.
Let $B=G/P$ be the flag manifold associated to $G$ as in Lemma~\ref{L:semi-simple}.
Then $(B,\nu)$ is a $G$-boundary as well as a $\Gamma$-boundary 
(Proposition~\ref{prop:boundary-lattice}).
Let $W=W_{G,B}$ be the $k$-relative Weyl group, as in Proposition~\ref{P:W4ss}.
Amenability of the $\Gamma$-action on $B$ yields 
a $\Gamma$-equivariant map $\Phi:B\to \Prob(\Circle)$.
Let $S(\Phi)<W$ be the subgroup corresponding to $\Phi$ under the Galois correspondence.
By Theorem~\ref{T:cyclic}, the action of $W$ on $W/S(\Phi)$ factors through 
a non-trivial finite cyclic quotient group, while Proposition~\ref{P:key4simple}
asserts that this action is faithful. 
This contradicts the fact that the $k$-relative Weyl group is never cyclic
for $\mathbf{G}$ of higher $k$ rank.
\end{proof}

\begin{proof}[Sketch of the proof of Theorem \ref{T:A2}]
The argument is essentially the same as for higher rank lattices.
Start from an $\tilde{A}_{2}$ group $\Gamma$ and a homomorphism 
$\ro:\Gamma\to\Homeo_{+}(\Circle)$.
If $\ro$ is elementary, then $\ro(\Gamma)$ has a finite orbit, because 
$\Gamma$ has a finite Abelianization, implied by property (T)
(\cite{CartwrightMlotkowskiSteger:T:94}).
We shall now reach a contradiction, assuming $\ro$ is non-elementary.
Corollary~\ref{C:no-inv-meas}, up to a semi-conjugacy we may assume that 
$\ro$ is minimal and strongly proximal.

Now let $(B,\nu)$ be a $\Gamma$-boundary as in Proposition~\ref{P:A2bnd},
and $W\cong S_{3}$ a subgroup of $W_{G,B}$.
Amenability of the $\Gamma$-action on $B$ yields 
a $\Gamma$-equivariant map $\Phi:B\to \Prob(\Circle)$.
Let $S(\Phi)<W$ be the subgroup corresponding to $\Phi$ under the Galois correspondence.
Then Theorem~\ref{T:cyclic} implies that the action of $W$ on $W/S(\Phi)$ factors through 
a non-trivial finite cyclic quotient group, while Proposition~\ref{P:key4A2}
asserts that this action is faithful. The fact that $S_{3}$ is not cyclic leads to the contradiction.
\end{proof}

\section{Commensurator superrigidity}
\label{S:commensurator}
\begin{proof}[Proof of Theorem~\ref{T:commensurator}]
Recall the notations: $\Gamma<G$ a lattice, $\Gamma<\Lambda<{\rm Commen}_G(\Gamma)$
is a dense subgroup in $G$ and $\ro:\Lambda\to \Homeo_+(\Circle)$ is an action 
without invariant measures. 
Up to a semi-conjugacy we may assume that $\ro(\Lambda)$ is 
minimal and strongly proximal (Corollary~\ref{C:no-inv-meas}).

The assumption $\Lambda<{\rm Commen}_G(\Gamma)$ means that for each $\lambda\in\Lambda$ the groups
\[
	\Gamma^{\prime}_{\lambda}=\lambda^{-1}\Gamma\lambda \cap \Gamma
	\qquad\text{and}\qquad
	\Gamma^{\pprime}_{\lambda}=\Gamma \cap \lambda\Gamma\lambda^{-1} 
\]
have finite index in $\Gamma$ and $c_{\lambda}:g\mapsto \lambda g \lambda^{-1}$
is an isomorphism between them:
\begin{equation}
\label{e:lambda-conj}
	c_{\lambda}:\Gamma^{\prime}_{\lambda}	\overto{\cong} \Gamma^{\pprime}_{\lambda}.
\end{equation}

\begin{prop}
\label{P:Lambda2Gamma}
If $\ro: \Lambda\to\Homeo_+(\Circle)$ is minimal and 
strongly proximal, and $\ro(\Gamma)$ has only infinite orbits,
then the restriction $\ro|_{\Gamma}$ to $\Gamma\to\Homeo_{+}(\Circle)$
is also minimal and strongly proximal.
\end{prop}
\begin{proof}

We first claim that minimality of the $\Lambda$-action implies minimality of the
$\Gamma$-action. 
First consider an arbitrary subgroup $\Gamma_{0} <\Gamma$ of finite index.
Then $\Gamma_{0}$ and $\Gamma$ have the same minimal set $K$.
To see this consider $\Gamma_{1}=\cap_{\gamma\in\Gamma} \gamma \Gamma_{0}\gamma^{-1}$, which is a normal subgroup of finite index in $\Gamma$.
The groups $\Gamma_{1}<\Gamma_{0}<\Gamma$ have only infinite orbits.
Let $K_{1}$, $K_{0}$, and $K$ denote \emph{the} minimal sets for the actions of these groups
on the circle (Proposition~\ref{P:minimal-sets}).
Since $K$ is invariant under $\Gamma_{0}$, it contains \emph{a} $\Gamma_{0}$-minimal
set, which has to be $K_{0}$ by uniqueness. Similarly for $\Gamma_{1}<\Gamma_{0}$,
giving the inclusions $K_{1}\subset K_{0}\subset K$.
On the other hand $\Gamma_{1}$ is normal in $\Gamma$. 
Therefore $K_{1}$ is a $\Gamma$-invariant compact set, which implies
$K\subset K_{1}$. This means $K_{1}=K_{0}=K$ as claimed.

Any $\lambda\in\Lambda$ gives a conjugation  
$c_{\lambda}$ as in (\ref{e:lambda-conj}) between finite index subgroups of $\Gamma$. 
Thus $\lambda$ maps objects uniquely associated to $\Gamma_{\lambda}^{\pprime}$
to similar objects associated to $\Gamma_{\lambda}^{\pprime}$.
Applying this to minimal sets we get $\ro(\lambda)K=K$ for every $\lambda\in\Lambda$.
As $\Lambda$ acts minimally on the circle, we deduce $K=\Circle$, that is
the minimality of the $\Gamma$-action.

Next observe that for a minimal action of some group on the circle,
the alternatives of Theorem~\ref{T:Margulis-Ghys} can be phrased in terms of the
\emph{centralizer} of the action: it is either (1) conjugate to the full group of
rotations in case of an equicontinuous action, or (2) is a finite cyclic group.
For minimal group actions strongly proximality is equivalent 
to the triviality of the centralizer $C$.

Given a finite index subgroup $\Gamma_{0}$ in the lattice $\Gamma$,
and a normal subgroup $\Gamma_{1}$ as before, let  $C$, $C_{0}$ and $C_{1}$ 
denote the centralizers in $\Homeo_{+}(\Circle)$ of $\ro(\Gamma)$, 
$\ro(\Gamma_{0})$ and $\ro(\Gamma_{1})$ respectively. 
Clearly $C<C_{0}<C_{1}$. 
Normality of  $\Gamma_{1}$ in $\Gamma$ implies that $C_{1}$ is 
normalized by $\Gamma$.
Now observe that for both the full group of rotations, and for any
finite cyclic group $C$, the normalizer of $C$ in $\Homeo_{+}(\Circle)$
coincides with its centralizer.
Thus $C_{1}<C$.

Similarly to the minimality argument, for $\lambda\in\Lambda$ we get 
that $\ro(\lambda)$ conjugates the centralizer of $\Gamma_{\lambda}^{\prime}$
to that of $\Gamma_{\lambda}^{\pprime}$. 
But by the above argument both are $C$ -- the centralizer of $\Gamma$.
Hence $C$ is normalized by $\ro(\Lambda)$. 
It is then centralized by $\ro(\Lambda)$. 
The latter is minimal and strongly proximal, so $C=\{1\}$. 
This proves that $\Gamma$ is also minimal and strongly proximal.
\end{proof}

\medskip

Let $(B,\nu)$ be a $G$-boundary, and $\Phi:B\to\Prob(\Circle)$ be the $\Gamma$-equivariant 
boundary map. 
By Theorem~\ref{T:reduction2atoms} (1) and (2) such $\Phi$ is unique and has the form
\[
	\Phi(x)=\frac{1}{k} \sum_{a\in A_{x}}\delta_{a}
\]
for a measurable map $B\to\Circle_{k}$, $x\mapsto A_{x}$. 
In what follows we are suppressing $\ro$ from the action notation.
\begin{lem}
\label{L:equiv}
For any fixed $\lambda\in\Lambda$ let $\Phi_\lambda:B\to\Prob(\Circle)$ be defined by
\[
	\Phi_\lambda(x)=\frac{1}{n} \sum_{i=1}^n g_i^{-1}\lambda^{-1}\Phi(\lambda g_i.x)
\]
where $g_1,\dots,g_n\in\Gamma$ be the coset representatives for 
$\Gamma^\prime\backslash\Gamma$.
Then $\Phi_{\lambda}$ is $\Gamma$-equivarinat. 
\end{lem}
\begin{proof}
Take $\gamma\in\Gamma$. It permutes the cosets $\Gamma^\prime g_i$ by
$\Gamma^\prime g_i\gamma=\Gamma^\prime g_{\pi(i)}$
where $\pi$ is a permutation of $\{1,\dots,n\}$.
In particular there exist $\gamma_i^\prime\in\Gamma^\prime$ so that
\[
	g_i \gamma=\gamma^\prime_i g_{\pi(i)}.
\]
Denote $\gamma_i^\pprime=\lambda\gamma_i^\prime\lambda^{-1}\in
\Gamma^{\pprime}_{\lambda}<\Gamma$. 
\begin{eqnarray*}
	\Phi_\lambda(\gamma.x)&=& 
		\frac{1}{n} \sum_{i=1}^n g_i^{-1}\lambda^{-1}\Phi(\lambda g_i\gamma.x)
		=\frac{1}{n} \sum_{i=1}^n g_i^{-1}\lambda^{-1}\Phi(\lambda \gamma_i^\prime g_{\pi(i)}.x)\\
		&=&\frac{1}{n} \sum_{i=1}^n g_i^{-1}\lambda^{-1}\Phi(\gamma_i^\pprime\lambda g_{\pi(i)}.x)
		=\frac{1}{n} \sum_{i=1}^n g_i^{-1}\lambda^{-1}\gamma_i^\pprime\Phi(\lambda g_{\pi(i)}.x)\\
		&=&\frac{1}{n} \sum_{i=1}^n \gamma g_{\pi(i)}^{-1}\lambda^{-1}\Phi(\lambda g_{\pi(i)}.x)
		=\gamma \Phi_\lambda(x).
\end{eqnarray*}
\end{proof}
Recall that Theorem~\ref{T:reduction2atoms} asserts that 
\[
	\Phi(x)=\frac{1}{k} \sum_{a\in A_{x}}\delta_{a}
\]
is \emph{the unique} $\Gamma$-equivariant map $B\to\Prob(\Circle)$
(here $x\mapsto A_{x}$ is a measurable $\Gamma$-equivariant map $B\to\Circle_{k}$).
Hence $\Phi_{\lambda}=\Phi$. Write $\Phi_{\lambda}$ as
\[
	\Phi_{\lambda}(x)=\frac{1}{n}\sum_{i=1}^{n} \Phi_{\lambda,i}(x),\qquad
	\Phi_{\lambda,i}(x)=g_{i}^{-1}\lambda^{-1}\Phi(\lambda g_{i}.x).
\] 
Each of $\Phi_{\lambda,i}(x)$ is supported on a $k$-point set, and their average is also
supported on a $k$-point set $A_{x}$.
As all $\Phi_{\lambda,i}(x)$ and $\Phi(x)$ are equidistributed on their support we get
\[
	\Phi_{\lambda,i}(x)=\Phi(x)\qquad \forall i=1,2,\dots,n
\]
for $\nu$-a.e. $x\in B$. For any $i$ we have
\[
	\lambda^{-1} \Phi(\lambda.(g_{i}.x))=g_{i} \Phi(x)=\Phi(g_{i}.x)
\]
because $g_{i}\in\Gamma$. As $g_{i}\nu\sim\nu$ we get an $\nu$-a.e. identity 
\[
	\lambda^{-1}\Phi(\lambda.y)=\Phi(y)\qquad\text{i.e.}\qquad \Phi(\lambda. y)=\lambda \Phi(y).
\]
We have proved that $\Phi:B\to\Prob(\Circle)$ is $\Lambda$-equivariant.

\medskip

We can now apply Lemma~\ref{L:extending} to
deduce that $\ro$ extends (uniquely) to a continuous homomorphism
\[
	\tilde\ro:G\to\Homeo_{+}(\Circle)
\]
and will denote by $L=\bar\ro(G)$ its image. 
Recall that this is a non-discrete \lcsc subgroup acting minimally and strongly proximally
on the circle, and hence by Theorem~\ref{T:lcsubgroups} $L$ is conjugate to $\PSL_{2}(\RR)$.

If $G$ is topologically simple then the epimorphism $p:G\to\PSL_{2}(\RR)$ is an isomorphism.
It then follows from Margulis commensurator superrigidity that the lattice $\Gamma$ has to
be arithmetic, since it has a dense commensurator $\Lambda$ in $\PSL_{2}(\RR)$.
This completes the proof of the theorem.
\end{proof}

\section{Cocycle versions of the results}
\label{S:cocycles}
%

Let $(X,\mu)$ be a standard probability space, $G$ a \lcsc group or a discrete group
acting ergodically by measure preserving transformations on $(X,\mu)$,
and let $\ro:G\times X\to \Homeo_{+}(\Circle)$ be a measurable cocycle.
In what follows we shall need to consider spaces of measurable functions, or sections,
from $X$ to $Y$ where $Y$ represents a topological space such as 
$\Circle$, $\Prob(\Circle)$, or $\Homeo_{+}(\Circle)$ etc.
More precisely we denote by $\Sect(X,Y)$ the set of equivalence classes
of measurable functions $\sigma:X\to Y$, where $\sigma\sim\sigma'$ iff
$\mu\setdef{x}{\sigma(x)\neq\sigma'(x)}=0$. 
All $\Sect(X,Y)$ can be endowed with a natural Borel structure, 
and some such spaces with an additional structure. 
For example $\Sect(X,\Prob(\Circle))$ can be viewed as a convex compact set 
with an affine $G$-action $a:G\to\Aff(\Sect(X,\Prob(\Circle)))$
\[
	(a(g) \sigma)(g.x)=\ro(g,x)_{*}\sigma(x)\qquad
	(\sigma\in\Sect(X,\Prob(\Circle)))
\]
so that the $G$-action is continuous (cf. \cite{Furstenberg:afterMargulisZimmer:81},
and Appendix~\ref{S:Appendix} below)\footnote{ 
The compact metric topology on $\Sect(X,\Circle)$, or $\Sect(X,\Prob(\Circle))$,
corresponds to convergence in measure (also in $L^{1}(X,\mu)$) of measurable functions
$X\to \Circle$ or $X\to\Prob(\Circle)$.}.

\begin{thm}[Cocycle version of Theorem~\ref{T:Margulis-Ghys}, \ref{C:no-inv-meas}]
\label{T:cocycle-min-strprox}
Let $\ro:G\times X \to\Homeo_{+}(\Circle)$ be a non elementary cocycle.
Then $\ro$ is semi-conjugate to a measurable cocycle 
$\ro':G\times X \to\Homeo_{+}(\Circle)$ which is minimal 
and strongly proximal.
\end{thm}
The precise definitions of minimality and strong proximality for cocycles
and the proof of the above theorem are given in the Appendix below 
(see Definitions~\ref{D:cocycle-minstr} and the discussion afterwards).
This result implies the following cocycle analogue of Theorem~\ref{T:reduction2atoms};
the proof of which is a straightforward cocycle adaptation of the proof  of 
Theorem~\ref{T:reduction2atoms}.
\begin{thm}
\label{T:Phi-coc}
Let $\ro:G\times X\to\Homeo_{+}(\Circle)$ be a minimal and strongly proximal 
measurable cocycle.
Let $(B,\nu)$ be a $G$-boundary, and $\Phi:B \to \Sect(X,\Prob(\Circle))$
a measurable $G$-equivariant map.
Then
\begin{enumerate}
\item
	Such $\Phi$ is unique (as a map defined up to $\nu$-null sets).
\item
	For $\nu\times\mu$-a.e. $(b,x)$ the measure $\Phi(b)(x)$ is atomic, 
	supported and equidistributed on a set $A_{b,x}\subset \Circle$ of $k$ points.
\item
	For a.e. $(b_{1},b_{2},x)\in B\times B\times X$ 
	the two-point sets $A_{b_{1},x}$ and $A_{b_{2},x}$ 
	are "unlinked" that is they belong to complimentary arcs of the circle.
\item
	In fact, $k=1$ or $k=2$.
\end{enumerate}
\end{thm}
Recall the definitions and the construction of an abstract Weyl group in \ref{SS:Weyl}. 
In an analogy with Theorem~\ref{T:cyclic} we have the following proposition.
\begin{prop}
\label{P:cocycle-cyclic}
Let $\ro:G\times X\to\Homeo_{+}(\Circle)$ be a minimal and strongly proximal 
measurable cocycle.
Let $(B,\nu)$ be a $G$-boundary, and $\Phi:B \to \Sect(X,\Prob(\Circle))$ be
a measurable $G$-equivariant map.
Let $W$ be a finite subgroup of $W_{G,B}$, containing $w_0$.
Then $W$ acts on $W/S(\Phi)$ through a non-trivial finite cyclic quotient group.
\end{prop}
\begin{proof}
By Theorem~\ref{T:Phi-coc}, the image of $\Phi$ is in $\Sect(X,\Circle_k)$.
Consider the map 
\[ 
	\vec{\Psi}:B\times B \overto{} \prod_{W/S(\Phi)} \Sect(X,\Circle_{k}) 
		\cong \Sect(X,(\Circle_{k})^{m})
\]
given by Construction \ref{cons:vec-map}, where we denote $|W/S(\Phi)|=m$.
The associated composition
\[ 
	B \times B \times X \overto{} (\Circle_{k})^{m} \overto{} 
	(\Circle_{k})^{m}/\Homeo_{+}(\Circle) 
\]
must be essentailly constant, by the ergodicity of $B\times B \times X$ and Lemma~\ref{l:tame}.
Hence it lies in a single $\Homeo_{+}(\Circle)$ orbit $\mathcal{O}\subset (\Circle_{k})^{m}$.
The corresponding $W$-equivariant map
\[
	B\times B \times X\overto{} \mathcal{O} 
\]
pushes the measure class on $B\times B\times X$ to a $W$-invariant measure class on 
$\mathcal{O}$.
By Lemma~\ref{l:aut(orbit)}, we get that the action of $W$ on $\mathcal{O}$ 
factors through a finite cyclic quotient.
The proof is complete, as the $W$-action on $\mathcal{O}$ is given by a coordinate permutation,
obtains from its action on $W/S(\Phi)$.
We are just left to remark that since $\ro$ is non-elementary,
$\Phi$ is a non-trivial factor map, hence by Lemma~\ref{l:long-element},
$w_0 \notin S(\Phi)$, so the $W$ action on $W/S(\Phi)$ is not trivial.
\end{proof}

\medskip

\subsection{Proofs of Theorems \ref{T:WZ} 
and \ref{T:cocycle-A2}} 

The proofs of these cocycle results differ from those of 
Theorems~\ref{T:higherrank} and  \ref{T:A2} only in replacing
the results about homomorphisms by their cocycle analogues.
More specifically the argument is as follows.

Let $G$ denote a higher rank simple group over a local field, 
a lattice there or a (discrete) $\tilde{A}_{2}$ group, and 
let $\ro:G\times X\to\Homeo_{+}(\Circle)$ be a measurable cocycle 
that we assume to be non-elementary. 
Upon applying a semi-conjugacy we may assume that $\ro$ is 
minimal and strongly proximal (Theorem~\ref{T:cocycle-min-strprox}).
Let $(B,\nu)$ be an appropriate $G$-boundary: as in 
Lemma~\ref{L:semi-simple} for semi-simple $G$, and as in 
Proposition~\ref{P:A2bnd} for the $\tilde{A}_{2}$-groups.
Let $W=W_{G,B}$ be the $k$-relative Weyl group as in Propositions~\ref{P:W4ss},
or $W\cong S_{3}$ as in \ref{P:W4A2}, respectively.
By amenability of the $G$-action on $(B,\nu)$ there exists a $G$-equivariant
map $\Phi:B\to\Sect(X,\Prob(\Circle))$.
Theorem~\ref{P:A2bnd}  asserts that it can be viewed as
$\Phi:B\to\Sect(X,\Circle_{k})$ for $k=1$ or $k=2$.
In any case let $S(\Phi)<W$ denote the subgroup corresponding to this quotient of $B$
in the Galois correspondence between $\mathcal{SG}(W)$ -- subgroups of $W$, and 
$\mathcal{Q}(B)$ -- the quotients of $B$. 
Proposition~\ref{P:cocycle-cyclic} asserts that the $W$-action on $W/S(\Phi)$
factors through a non-trivial cyclic group. 
Complexity of $W$ precludes this from happening.
Indeed for a simple algebraic case, Proposition~\ref{P:key4simple}
states that such an action is always faithful, and this leads a contradiction
as the $k$-relative Weyl group of a simple $k$-group $\mathbf{G}$ with 
$\text{rk}_{k}(\mathbf{G})\ge 2$ is never cyclic.
In the $\tilde{A}_{2}$ case the contradiction is that $S_{3}$ is not cyclic,
while $S_{3}$ acts faithfully on $W/S(\Phi)$ (Proposition~\ref{P:key4A2}).
 
\medskip
 
\subsection{Proof of Theorem \ref{T:cocycle-product}}
Let $G=G_{1}\times \cdots \times G_{n}$ product of some \lcsc groups,
acting measurably on a probability space $(X,\mu)$, so that each $G_{i}$ acts ergodically.
Let $\ro:G\times X\to\Homeo_{+}(\Circle)$ be a non-elementary measurable cocycle.
Note that not all of $G_{i}$ are amenable, for otherwise $G$ would be amenable and any cocycle 
would be elementary. 
Up to reordering, we may assume that for some $1\le m\le n$ the factors
$G_{1},\dots,G_{m}$ are non-amenable while $G_{i}$ are amenable for $m<i\le n$.

Let  $B=B_{1}\times \cdots \times B_{m}$ be the product of the non-trivial boundaries 
as in Proposition~\ref{P:products-of-boundaries}, and $W=(\ZZ/2\ZZ)^{m}<W_{G,B}$ 
be a subgroup of the Weyl group as in \ref{O:Weyl4products}.
By Theorem~\ref{T:cocycle-min-strprox}, up to a semi-conjugacy $\ro$ can be assumed
to be minimal and strongly proximal.
Let $\Phi:B\to \Sect(X,\Prob(\Circle))$ be the corresponding boundary map,
the existence of which is guaranteed by the amenability of $B$.
By Theorem~\ref{T:Phi-coc}, the image of $\Phi$ is in $\Sect(X,\Circle_k)$
for $k=1$ or $k=2$.
Proposition \ref{P:cocycle-cyclic} shows that the action of $W\cong (\ZZ/2\ZZ)^{m}$
on $W/S(\Phi)$ factors through a non-trivial finite cyclic group.
By Lemma~\ref{L:WPhi4products}, $S(\Phi)=W_I$ 
for $I=\{1,\ldots,m\}-\{i\}$ for some $i\in\{1,\ldots,m\}$, and up to reordering we may assume 
that $i=1$.
Hence $\Phi$ factors through $\Phi_{1}:B_{1}\to \Sect(X,\Circle_{k})$.
We denote $H=G_{2}\times\dots\times G_{n}$, and observe that $H$
acts trivially on $B_{1}$, 
but ergodically on $(X,\mu)$.  

\medskip

Recall that at a similar stage in the proof of Theorem~\ref{T:product}, we were
able to continuously \emph{extend the given} action of a dense subgroup 
in $G_{1}$, namely $\pi_{1}(\Gamma)$, to $\ro_{0}:G_{1}\to\Homeo_{+}(\Circle)$.
The challenge in the present setting is that there is no circle action to start with! 
We shall construct an ``abstract'' circle and a $G_{1}$-action on it,
and will prove that the cocycle $\ro$ is conjugate to this homomorphism.
Theorem~\ref{T:cocycle-product} will be proved in four steps.

\medskip

{\bf Step 1:} \emph{There exists a probability measure $\eta$ on the compact metric space
$\Sect(X,\Circle)$ so that }
\begin{itemize}
\item[(i)]
	\emph{for all $g\in G$ we have $a(g)_{*}\eta\sim\eta$},
\item[(ii)]
	\emph{for all $g\in H$ we have $a(g)(\sigma)=\sigma$ for $\eta$-a.e. 
	$\sigma\in\Sect(X,\Circle)$.}
\end{itemize} 
Recall that we have a measurable map $\Phi:B\to\Sect(X,\Circle_{k})$ for $k=1$ 
or $k=2$, which factors through $B_{1}$:
\[
	\Phi:B\overto{}B_{1}\overto{\Phi_{1}}\Sect(X,\Circle_{k})
\]
and such that for each $g\in G$:
$\Phi(g.b)=\Phi_{1}(\pi_{1}(g).b_{1})=g.\Phi_{1}(b_{1})=g.\Phi(b)$ 
for $\nu$-a.e. $b\in B$ and $\nu_{1}$-a.e. $b_{1}\in B_{1}$.
If $k=1$ then the pushforward measure
\[
	\eta=\Phi_{*}\nu=(\Phi_{1})_{*}\nu_{1}
\]
clearly satisfies (i) and (ii) above.

The case of $k=2$ is a little trickier.
Instead of $\Phi:B\to\Sect(X,\Circle_{2})$, consider the 
map $\Psi(b',b'')=(\Phi(b'),\Phi(b''))$
\[
	\Psi:B\times B\to \Sect(X,\Circle_{2}\times \Circle_{2})\cong
	\Sect(X,\Circle_{2})\times\Sect(X,\Circle_{2}).
\]
A.e. on $B\times B\times X$ the two point sets $A_{b',x}=\Phi(b')(x)$,  
$A_{b'',x}=\Phi(b'')(x)$ are disjoint and ``unlinked'' (Theorem \ref{T:Phi-coc} (3)). 
We claim that $\Psi$, viewed as a four point bundle over $B\times B\times X$, 
can be ``trivialized''. More precisely we can define measurable 
$G$-equivariant maps $\Psi_{i}:B\times B\to \Sect(X,\Circle)$, $i\in\{1,2,3,4\}$, 
so that
\begin{eqnarray*}
	\Phi(b')(x)&=&\{\Psi_{1}(b',b'',x), \Psi_{2}(b',b'',x)\},\\ 
	\Phi(b'')(x)&=&\{\Psi_{3}(b',b'',x), \Psi_{4}(b',b'',x)\}.
\end{eqnarray*}
Indeed,  a four point set $P\subset \Circle$ with a a given partition $P=P'\cup P''$ into two 
``unlinked'' pairs $P'$ and $P''$ of points, can be written in a unique way as 
$P=\{p_{1},p_{2},p_{3},p_{4}\}$ so that the cyclic order of the 
points would match the cyclic order $(1,2,3,4)$ of their indices, while the partition 
would be into $P'=\{p_{1},p_{2}\}$ and $P''=\{p_{3},p_{4}\}$.
Applying this argument to the four point set $P_{x}=A_{b',x}\cup A_{b'',x}$
we obtain the above maps $\Psi_{i}(b',b'',x)$ (measurability being clear from the construction).
We can now use, say $\Psi_{1}$, to pushforward the measure $\nu\times\nu$ 
\[
	\eta=(\Psi_{1})_{*}(\nu\times\nu).
\]
Step 1 is accomplished.

\medskip

{\bf Step 2:} \emph{The support $\Sigma=\supp(\eta)\subset\Sect(X,\Circle)$ 
has a $G$-invariant dense cyclic order, compatible with $\mu$-a.e. 
evaluation projection $p_{x}:\Sigma\to\Circle$.}

\medskip

Let $\omega:\Circle\times\Circle\times\Circle\to\{-1,0,1\}$ denote the given
orientation on the circle ($\omega(a,b,c)=1$ represents  
$(a,b,c)=(b,c,a)=(c,a,b)$, $\omega(a,b,c)=-1$ represents
the order $(a,c,b)=(c,b,a)=(b,a,c)$, and we set $\omega(a,b,c)=0$ if not all the points are distinct).
An abstract cyclic order on $\Sigma$ is a map $\Omega:\Sigma^{3}\to\{-1,0,1\}$
satisfying some obvious consistency relation (algebraically, the condition is that $\Omega$ 
is a cocycle: $\pa\Omega$ vanishes on $\Sigma^{4}$).
We shall define such an order by pulling back $\omega$ via $\mu$-a.e. 
evaluation projection $p_{x}:\Sigma\to\Circle$.

Observe that since $G$-act continuously on $\Sect(X,\Circle)$ preserving
the measure class of $\eta$, the support $\Sigma=\supp(\eta)$ is a $G$-invariant set,
on which $H$ acts trivially.
Hence for any triple $\alpha,\beta,\gamma\in\Sigma$ and for any $g\in G$
we have for $\mu$-a.e. $x\in X$\footnote{Hereafter we simplify the notation 
for the $G$-action on $\Sigma$ by writing $g\sigma$ instead of $a(g)(\sigma)$}:
\begin{eqnarray*}
	 \omega(g\alpha(g.x),g\beta(g.x),g\gamma(g.x))&=&
	 \omega(\ro(g,x)\alpha(x),\ro(g,x)\beta(x),\ro(g,x)\gamma(x))\\
	 &=&\omega(\alpha(x),\beta(x),\gamma(x)).
\end{eqnarray*}
In particular, the measurable function $\omega\circ p_{x}$ is $H$-invariant
on $\Sigma^{3}$, and hence is $\mu$-a.e. constant. This allows to define 
the cyclic order $\Omega$ on $\Sigma$.
The above $\mu$-a.e. relation also shows that $G$ preserves $\Omega$.
It follows from minimality and strong proximality of $\ro$ that $\Omega$ is dense: 
for any $\alpha\neq\beta$ and $\gamma$ in $\Sigma$
there exists $g\in G$ so that 
\[
	E=\setdef{x\in X}{ \omega(\alpha(x),g\gamma(x),\beta(x))=1}
\]
has $\mu(E)>0$ (see Lemma~\ref{L:cocycle-minstr}), and therefore
$\Omega(\alpha,g\gamma,\beta)=1$ (and $\mu(E)=1$).

\medskip

{\bf Step 3:} \emph{$\Sigma\cong\Circle$ which gives a $G$-action:
$\ro_{0}:G\overto{\pi_{1}}G_{1}\overto{}\Homeo_{+}(\Circle)$.}

\medskip

Indeed, fix some countable subset $\Sigma_{0}\subset\Sigma$ which is 
topologically dense in $\Sigma$ and inherits a dense cyclic order, 
and choose a similar subset $S_{0}\subset\Circle$
(for example $S_{0}=\QQ/\ZZ\subset\Circle=\RR/\ZZ$).
Then there exists an order isomorphism $h_{0}:\Sigma_{0}\to S_{0}$.
By topological density of the sets this order isomorphism 
uniquely extends to a homeomorphism $h:\Sigma\to \Circle$. 
The $G$-action on $\Sigma$ then defines the action $\ro_{0}:G\to\Homeo_{+}(\Circle)$
by
\[
	\ro_{0}(g)h(\sigma)=h(g\sigma)\qquad(\sigma\in\Sigma,\ g\in G).
\]
This $G$-action factors through $G_{1}$ because $H$ act trivially on 
$\Sigma$.

\medskip

{\bf Step 4:} \emph{There exist $f\in\Sect(X,\Homeo_{+}(\Circle))$ s.t. 
$\ro(g.x)=f_{g.x}\circ \ro_{0}(g)\circ f_{x}^{-1}$.}

\medskip

Note that since $\Sigma$ is $G$-invariant, while $\ro$ is minimal,
$p_{x}(\Sigma)=\Circle$ for $\mu$-a.e. $x$.
The set $\Sigma_{0}$ is countable, hence there is a full measure set 
$X_{0}\subset X$ so that $\Omega$ on $\Sigma_{0}$ agrees with 
$\omega$ on $S_{x}=p_{x}(\Sigma_{0})\subset\Circle$.
The assumption that $\Sigma_{0}$ is dense in $\Sigma$ implies that for
any $\epsilon>0$:
\[
	\mu\left(\setdef{x\in X_{0}}{\Circle\setminus S_{x}\ 
	\text{contains an }\epsilon\text{-arc} }\right)=0.
\]
Thus for a subset $X_{1}\subset X_{0 }$ of full measure, $S_{x}$ is dense 
in $\Circle$, for all $x\in X_{1}$. 
Let $f_{x}:S_{0}\to S_{x}$ be an order preserving bijection; it can be chosen
measurably on $X_{1}$. 
We continue to denote by $f_{x}:\Circle\to\Circle$ the unique continuous
extension of $f_{x}:S_{0}\to S_{x}$. 
Then $\{f_{x}\}_{x\in X_{1}}$ is a measurable family of 
homeomorphisms, such that for $\sigma\in\Sigma$ for $\mu$-a.e. $x\in X$:
$\sigma(x)=(f_{x}\circ h)(\sigma)$.
Thus for any $g\in G$ and any countable set $\{\sigma_{i}\}$ in $\Sigma$,
 we have for $\mu$-a.e. $x\in X$:
\[
	\ro_{0}(g) (p_{i}) = (f_{g.x}^{-1}\circ \ro(g,x)\circ f_{x}) (p_{i}),
	\qquad\text{where}\qquad p_{i}=h(\sigma_{i}).
\]
Choosing the sequence $\{\sigma_{i}\}$ to be dense in $\Sigma$ this 
gives the identity of homeomorphisms: for every $g$ for $\mu$-a.e. $x\in X$:
$
	\ro_{0}(g)=f_{g.x}^{-1}\circ \ro(g,x)\circ f_{x},
$
where $\ro_{0}:G\overto{\pi_{1}}G_{1}\overto{}\Homeo_{+}(\Circle)$.
This complets the proof of Theorem~\ref{T:cocycle-product}.

\section{Appendix}
\label{S:Appendix}
\subsection{Generalities on function spaces}

Let $(X,\mu)$ be a Lebesgue probability space and $M$ be a topological space.
We denote by $\Sect(X,M)$  the space of equivalence classes of measurable 
functions $X\to M$, defined up to measure zero.
We endow $\Sect(X,M)$ with the weakest topology such that for every $f\in C_c(M)$,
and $\lambda \in L^1(X,\mu)$, the function
\begin{equation}
\label{e:wtop}
	\Sect(X,M) \ni \sigma\ \mapsto\ \int_{X} \lambda(x) \cdot f(\sigma(x))\,d\mu
\end{equation}
is continuous. The following are straightforward but useful observations: 
\begin{lem}
 
\begin{enumerate}
\item 
	$\Sect(X,M)$ is functorial:  contravariant in $G$-actions on $X$, 
	and covariant in $M$.
\item 
	$\Sect(X,M_1) \times \Sect(X,M_2) \cong \Sect(X,M_1\times M_2)$.
\item 
	$\Sect(X_1,\Sect(X_2,M)) \cong \Sect(X_{1}\times X_2,M)$.
\item 
	If $M$ is a topological group then $\Sect(X,M)$ is a topological group
	(but even if $M$ is locally compact, $\Sect(X,M)$ need not be Polish).
\item 
	If $M$ is compact (and metrizable) then $\Sect(X,M)$ is compact (and metrizable).
\item 
	If $M$ is a compact convex subset of a topological vector space then
	$\Sect(X,M)$ can be seen as a compact convex subset of a topological vector space.
\end{enumerate}
\end{lem}

Let $M_1$ be a locally compact, second countable space, and $M_2$ be a Polish space.
The space $C(M_1,M_2)$, endowed with the compact-open topology is known to be Polish.
Another natural topology on $C(M_1,M_2)$ is the topology of pointwise convergence.
The following lemma is easy.

\begin{lem}
Let $M_1$ be a locally compact, second countable space, and $M_2$ be a Polish space.
The topology of pointwise convergence and the open compact topology on
$C(M_1,M_2)$ generate the same $\sigma$-algebra.
\end{lem}

\begin{rem}
In our special cases of interest, $\Homeo_+(\RR)$ and $\Homeo_+(\Circle)$,
the two topologies actually coincide (Dini's theorem).
\end{rem}

\begin{cor} \label{c:measurablility}
For an almost everywhere defined function $\sigma:X\to C(M_1,M_2)$,
the following are equivalent.
\begin{enumerate}
\item
	$\sigma\in \Sect(X,C(M_1,M_2))$.
\item 
	For any open $U\subset M_2$ and $p\in M_1$,
	the set $\{x\in X \mid  \sigma(x)(p)\in U \}$ is measurable in $X$.
\item
	For every $p\in M_1$, the section $x \mapsto \sigma(x)(p)$ is in $\Sect(X,M_2)$.
\item
	For every $\tau\in \Sect(X,M_1)$, the section $x \mapsto \sigma(x)(\tau(x))$ is in $\Sect(X,M_2)$.
\end{enumerate}
\end{cor}

\begin{lem}
\label{L:cyclic-bundle}
Let $E$ be a measurable bundle over $X$, with countable fiber,
and assume that on (almost) every fiber there is a dense cyclic order, 
where the order relation is measurable.
Then there is a conull set $X_{0}\subset X$ and an isomorphism of bundles 
over $X_{0}$ between $E$ and $X_{0}\times \QQ/\ZZ$,
which is order preserving on almost every fiber.
\end{lem}
\begin{proof}
The proof of this statement is a straightforward modification of the standard
proof that every countable dense cyclic order is isomorphic to that of $\QQ/\ZZ$.
One just needs to replace terms as "choose a point such that..."
by "chosse a measurable section such that...".
\end{proof}

\subsection{Circle Bundles}
The cocycle versions of the notions of \emph{minimality} and \emph{strong proximality}
were introduced in Furstenberg's \cite{Furstenberg:afterMargulisZimmer:81}.
For circle bundles it is convenient to combine together these terms.
Denote by $2^M$ the compact metric space of all non-empty closed subsets of 
$M$ equipped with the Hausdorff metric.
\begin{defns}
\label{D:cocycle-minstr}
Let $\ro:G\times X \to \Homeo_+(\Circle)$ be a measurable cocycle over an ergodic
probability measure preserving $G$-action on $(X,\mu)$. We say that $\ro$ is:
\begin{itemize}
\item
	\emph{elementary} 
	if there exists a $G$-invariant section in $\Sect(X,\Prob(\Circle))$;
\item
	\emph{minimal} 
	if there are no non-trivial $G$-invariant sections in $\Sect(X,2^{\Circle})$;
\item
 	\emph{minimal and strongly proximal}
	if for any section of proper subsets, $\sigma\in \Sect(X,2^{\Circle})$
	and any $\tau\in\Sect(X,\Circle)$, 
	the $G$-orbit of $\sigma$ contains $j(\tau)$ in its closure,
	where $j(\tau)(x)=\{\tau(x)\}\in \Sect(X,2^{\Circle})$.
\end{itemize}
\end{defns}
To clarify the above definitions, here is an alternative and more concrete description,
essentially identical  to the one in \cite[Lemmas 3.2, 3.3]{Furstenberg:afterMargulisZimmer:81}.
\begin{lem}
\label{L:cocycle-minstr}
Let $\ro:G\times X\to \Homeo_{+}(\Circle)$ be a measurable cocycle over an ergodic
probability measure preserving $G$-action on $(X,\mu)$. Then
\begin{enumerate}
\item
	 $\ro$ is minimal iff given any $\sigma\in\Sect(X,\Circle)$, a non-empty open 
	 $U\subset\Circle$, and a measurable subset $E\subset X$ with $\mu(E)>0$, 
	 there exists $g\in G$ and $x\in E\cap g^{-1}E$ so that $\ro(g,x)\sigma(x)\in U$.
\item
	 $\ro$ is minimal and strongly proximal iff given any 
	 $K\in\Sect(X,2^{\Circle}\setminus\{\Circle\})$,
	 a non-empty open $U\subset\Circle$, and a measurable subset 
	 $E\subset X$ with $\mu(E)>0$, 
	 there exists $g\in G$ and $x\in E\cap g^{-1}E$ so that $\ro(g,x) K(x)\subset U$.
\end{enumerate}
\end{lem}
%

\begin{lem}[Reduction to Minimal cocycles]
\label{L:cocycle-min}
Every non-elementry cocycle $\ro:G\times X \to \Homeo_+(\Circle)$ 
is semi-conjugate to a non-elementry minimal cocycle.
\end{lem}
\begin{proof}
By Zorn's lemma $\Sect(X,2^{\Circle})$ contains a minimal $G$-invariant 
element $K$ (cf.~\cite[Lemma 3.2]{Furstenberg:afterMargulisZimmer:81}).
By the ergodicity of the $G$-action on $X$, the topological type of the closed sets 
$K(x)\subset \Circle$ is $\mu$-a.e. the same.
Let $I(x)\subset K(x)$ denote the set of all the isolated points of $K(x)$. 
Then $I=\{I_{x}\}_{x\in X}$ is a $G$-invariant relatively open subset of $K$.
Hence $I$ is either empty, or $I=K$. The latter would imply that $K(x)=I(x)$
is a finite set with a fixed number of elements. 
This would contradict the assumption that $\ro$ is non-elementary.
Hence $I(x)=\emptyset$, i.e. $K(x)$ is a perfect set a.e., which leaves two possibilities:
either $K(x)=\Circle$  or $K(x)$ is a Cantor set for $\mu$-a.e. $x\in X$.
We need to reduce the latter case to the former using a semi-conjugacy.

For a Cantor set $K(x)$ the compliment $\Circle\setminus K(x)$ is a countable union 
of disjoint open intervals, which have a dense cyclic order type. 
This defines a bundle as in Lemma~\ref{L:cyclic-bundle}, and up to a null set we can write
\[
	\Circle\setminus K(x)=\bigcup_{q\in\QQ/\ZZ} U_{q}(x)
\]
where $U_{q}(x)$ varies measurably in $x\in X$ for each $q\in\QQ/\ZZ$.
This defines a measurable family $f_{x}:\Circle\to\Circle$ of order preserving 
continuous surjections with $f_{x}(K(x))=\Circle$ and $f_{x}^{-1}(\{q\})=\ol{U}_{q}(x))$.
It also defines a $G$-action on $\Sect(X,\Circle)$ via a cocycle 
$\ro':G\times X\to \Homeo_{+}(\Circle)$ which acts minimally on $X\times \Circle$
and satisfies
\[
	\ro'(g,x)\circ f_{x}=f_{g.x}\circ\ro(g.x).
\]
Cocycle $\ro'$ is non-elementary, because any $\ro'$-invariant measure projecting to $\mu$
would have given rise to a $\ro$-invariant probability measure.
\end{proof}

\begin{proof}[Proof of Theorem~\ref{T:cocycle-min-strprox}]
Let $\ro:G\times X\to\Homeo_{+}(\Circle)$ be a non-elementary cocycle. 
Up to a semi-conjugacy we may assume $\ro$ to be minimal (Lemma~\ref{L:cocycle-min}). 
We shall now construct a finite quotient which will give rise to a minimal and 
\emph{strongly proximal} cocycle.
The main issue here is the construction of an element $\theta \in \Sect(X,\Homeo_{+}(\Circle))$
that commutes with the action of $G$ on $X\times \Circle$.
Note that while $G$ moves fibers over $X$ according to the ergodic $G$-action on $X$,
$\theta$ will leave each circle fiber invariant. 

The construction is a fiberwise imitation of the construction in the proof of 
Theorem~\ref{T:Margulis-Ghys}.
Let $\pi:\RR\to \RR/\ZZ=\Circle$ denote the projection map.
For every $(x,p) \in X\times \RR$, define 
\[
	 \tilde{\theta}(x,p)=\sup\{q\in\RR~|~\forall \epsilon>0~\exists g\in G,\ \text{s.t.}\ 
	 \lambda(\ro(g,x)(\pi([p,q]))< \epsilon\} 
\]
where $\lambda$ denotes the length on $\Circle=\RR/\ZZ$.
As in the proof of Theorem~\ref{T:Margulis-Ghys}, it is not hard to see that for 
$\mu$-a.e. $x\in X$: the map $\tilde{\theta}(x):\RR\to\RR$ is a continuous monotonically
strictly increasing map, commuting with the $\ZZ$-translation.
The map $x\to\tilde{\theta}(x,p)$ is measurable for any fixed $p\in \RR$,
so by Corollary~\ref{c:measurablility}, $\tilde{\theta}$ is in $\Sect(X,\Homeo+(\RR))$.
Since it commutes with the action of $\Sect(X,\ZZ)$, $\tilde{\theta}$ defines 
an element $\theta\in \Sect(X,\Homeo_{+}(\Circle))$.
By the construction $\theta$ commutes with the $G$-action on $X\times\Circle$.

For $(x,p)\in X\times\RR$ we have $p\le\tilde{\theta}(x,p)\le p+1$.
The extreme cases $\tilde\theta(x,p)=p$ and  $\tilde\theta(x,p)=p+1$,
either occur for $\mu$-a.e. $x\in X$ and all $p\in\RR$ or for no $p$.
It is also easy to see that the limit 
$\tau(x)=\lim_{n\to\infty}  \tilde{\theta}^n(x,p)/n$ exists and
is independent of $p$, and $\tau:X\to [0,1]$ is a measurable function.
Being $G$-invariant, this function is a.e. constant $\tau$.
There are several possibilities for the values for $\tau$:

\medskip

{\bf Case $\tau=1$.} This means that $\tilde{\theta}(x,p)=p+1$ for a.e. $x\in X$ and 
all $p\in\RR$. Recalling the definition of $\tilde{\theta}$ this means that $\ro$ 
is strongly proximal. 

\medskip

{\bf Case $\tau\in(0,1)\cap\QQ$.} 
Suppose that $\tau=p/q$ with $\gcd(p,q)=1$. Then for $\mu$-a.e. $x\in X$ the homeomorphism
$\theta(x)$ has a rotation number $\tau=p/q$ and so $\theta^{q}(x)$ has a fixed point.
Let $F(x)\subset\Circle$ denote the fixed point set of $\theta^{q}(x)$.
Then $F\in\Sect(X,2^{\Circle})$ is $G$-invariant, and by minimality $F(x)=\Circle$, 
i.e. $\theta^{q}(x)$ is the identity for a.e. $x\in X$.
In this case the $q$-to-$1$ continuous map $f_{x}:\Circle\to \Circle/\theta(x)$ provides 
a semi-conjugacy of $\ro$ to a cocycle $\ro'$ with $\tau=1$, i.e. to a strongly proximal
(and minimal) cocycle, as claimed.

\medskip

{\bf Case $\tau\in (0,1)\setminus\QQ$.} 
In this case, $\theta(x)$ has an irrational rotation number and therefore is (semi)-conjugate
to an irrational rotation using $\tilde{h}:X\to C(\RR,\RR)$ defined by
$h(x)(p)=\sup\{\tilde{\theta}(x)^n(p)-n\tau\}$, which descends  to $h:X\to C(\Circle,\Circle)$
In fact, minimality of the $G$-action implies that $h(x)\in \Homeo(\Circle)$,
so this is a conjugation rather then semi-conjugacy.
This leads to a contradiction because the unique $\theta(x)$-invariant probability measure
would give rise to a $G$-invariant element in $\Sect(X,\Prob(\Circle))$, contrary
to the assumption that $\ro$ is non-elementary.

\medskip

{\bf Case $\tau=0$.} 
Here the argument differs slightly from the one used in Theorem~\ref{T:Margulis-Ghys}.
In this case $\theta$ has a fixed points at almost every fiber.
By the minimality of $\ro$, $\theta$ is trivial.
Let $\delta$ denote some fixed metric on the circle $\Circle$.
Consider the following family of associated pseudo-metrics:
\[ 
	d_x=\inf_{g\in G} \ \delta\circ \ro(g,x).
\]
Obviously, $d_{g.x}=d_{x}\circ\ro(g,x)$.
The assumption $\tau=0$ means that for $\mu$-a.e. $x\in X$ we have
$\tilde{\theta}(x,p)=p$ for all $p$, which implies 
$d_{x}(p,q)>0$ for any $p\neq q$.
So $d_{x}$ is measurable family of actual metrics on $\Circle$.
This allows us to define the ``integrated'' metric $D$ on $\Sect(X,\Circle)$:
\[
	D(\sigma,\tau)=\int_{X}d_{x}(\sigma(x),\tau(x))\,d\mu(x).
\]
This metric is preserved by the $G$-action.
For $\mu$-a.e. $x$ the metric $d_{x}$ defines the usual topology on the circle $\Circle$.
Thus $D$ defines the same (compact metrizable) topology  on $\Sect(X,\Circle)$
as the weak topology introduced at the beginning of this section in (\ref{e:wtop});
because this topology describes convergence in measure (also in $L^{1}(X,\mu)$) 
of measurable functions $X\to\Circle$, a notion which coincides for $d_{x}$ and $\delta$.

Recall that the group of isometries of a metric compact space is a compact group
(follows from Arzela-Ascoli). So the closure $K$ of $G$ in the group of isometries
$\Isom(\Sect(X,\Circle),D)$ is compact. 
The pushforward of the Haar measure on $K$ to any $K$-orbit on 
$\Sect(X,\Circle)$ defines a $G$-invariant probability measure $M$
on $\Sect(X,\Circle)$.

Finally, observe that a $G$-invariant measure $M\in\Prob(\Sect(X,\Circle))$
defines a $G$-invariant $m\in \Sect(X,\Prob(\Circle))$, by
\[
	m_{x}=\int_{\Sect(X,\Circle)} \delta_{\sigma(x)}\,dM(\sigma)
\]
(Fubini guarantees that this definition is correct).
Existence of a $G$-invariant element $m\in\Sect(X,\Prob(\Circle))$ contradicts
the assumption that $\ro$ is non-elementary. 

\end{proof}

\def\cprime{$'$}
\providecommand{\bysame}{\leavevmode\hbox to3em{\hrulefill}\thinspace}
\providecommand{\href}[2]{#2}

\end{document}